\documentclass[10pt]{article}
\usepackage{fullpage,epsfig}
\usepackage{amsmath}
\usepackage{amsthm}
\usepackage{amssymb}
\usepackage{latexsym}
\numberwithin{equation}{section}
\newtheorem{example}{Example}[section]
\newtheorem{theorem}[example]{Theorem}
\newtheorem{proposition}[example]{Proposition}

\newtheorem{lemma}[example]{Lemma}

\newtheorem{remark}[example]{Remark}

\newcommand{\R}{{\mathbb R}}
\newcommand{\Grad}{\mathcal{GY}^{+\infty,-\infty}(\O;\R^{n\times n})}

\newcommand{\inv}{{\rm inv}}
\newcommand{\Rm}{{\R^{n\times n}}}
\newcommand{\N}{{\mathbb N}}
\renewcommand{\det}{{\rm det}\, }

\newcommand{\reg}{{\R^{n\times n}_{\rm inv}}}
\newcommand{\regpos}{{\R^{n\times n}_{\rm inv+}}}

\newcommand{\be}{\begin{eqnarray}}

\newcommand{\ee}{\end{eqnarray}}
\renewcommand{\r}{\varrho}
\renewcommand{\d}{{\rm d}}

\newcommand{\dm}{{\cal DM}^p_{\cal R}(\O;\R^{n\times n})}

\newcommand{\ra}{\right\rangle}

\newcommand{\la}{\left\langle}
\newcommand{\md}{{\rm d}}
\newcommand{\nd}{\nu_x(\d s)}

\renewcommand{\O}{\Omega}

\newcommand{\C}{{C_{0,{\rm inv}}(\R^{n\times n})}}

\newcommand{\A}[1]{\langle#1\rangle}
\newcommand{\Rn}{\R^{n}}

\newcommand{\wto}{\rightharpoonup}
\newcommand{\wtostar}{\stackrel{*}{\wto}}

\newcommand{\rca}{{\rm rca}}

\newcommand{\Rrho}{{R_\varrho^{n \times n}}}
\newcommand{\RrhoPlus}{{R_{\varrho+}^{n \times n}}}

\begin{document}
\begin{sloppypar}

\begin{center}
{\huge \bf Young measures supported on invertible matrices}\\
\vspace{3mm}
Barbora Bene\v{s}ov\'{a}$^{1,2}$, Martin Kru\v{z}\'{\i}k$^{3,4}$, Gabriel Path\'{o}$^{2,4}$ \\
\vspace{2mm}
$^1$ Institute of Thermomechanics of the ASCR, Dolej\v{s}kova 5, CZ-182 08 Praha 8, Czech Republic\\
$^2$ Faculty of Mathematics and Physics, Charles University, Sokolovsk\'{a} 83, CZ-186 75 Praha 8, Czech Republic\\
$^3$Institute of Information Theory and Automation of the ASCR, Pod vod\'{a}renskou
v\v{e}\v{z}\'{\i}~4, CZ-182~08~Praha~8, Czech Republic \\
$^4$ Faculty of Civil Engineering, Czech Technical
University, Th\'{a}kurova 7, CZ-166~ 29~Praha~6, Czech Republic

\bigskip

\bigskip
\bigskip
\begin{abstract}
Motivated by variational problems in nonlinear elasticity depending on the deformation gradient and its inverse,  we completely and explicitly describe Young measures generated by matrix-valued mappings   $\{Y_k\}_{k\in\N} \subset L^p(\O;\R^{n\times n})$, $\O\subset\R^n$, such that $\{Y_k^{-1}\}_{k\in\N} \subset L^p(\O;\R^{n\times n})$ is bounded, too. Moreover, the constraint $\det Y_k>0$ can be easily included and is reflected in a condition on the support of the measure.  This condition  typically occurs  in problems of nonlinear-elasticity theory  for hyperelastic materials  if $Y:=\nabla y$ for $y\in W^{1,p}(\O;\R^n)$.  Then we fully characterize the set of Young measures generated by gradients of a uniformly bounded  sequence in $W^{1,\infty}(\O;\R^n)$ where the inverted gradients are also bounded in $L^\infty\O;\R^{n\times n})$. This extends the original results due to D.~Kinderlehrer and P.~Pedregal \cite{k-p1}. 
\end{abstract}
\end{center}

\section{Introduction}
In this paper we investigate a new tool to study minimization problems for integral functionals defined over matrix-valued mappings that take values \emph{only} in the set of invertible matrices. Typical examples are found, e.g., in non-linear elasticity where static equilibria are  minimizers of the elastic energy, i.e., one is led to solve 
\be\label{motivation}
\mbox{ minimize }J(y):=\int_\O W(\nabla y(x))\,\md x\ ,
\ee
where $\O\subset\R^n$  denotes the reference configuration of the material, $y\in W^{1,p}(\O;\R^n)$ is the deformation, $1 < p \le +\infty$, $y=y_0$ on $\partial\O$,  and $W \colon \R^{n\times n}\to\R$ is the stored energy density, i.e., the potential of the first Piola-Kirchhoff stress tensor. Usually in elasticity, one demands  at least that $\det\nabla y\neq 0$ to assure the local invertibility of $\nabla y$ or even that $\det\nabla y>0$ in order to preserve orientation of  $y$.

If $W$ is polyconvex, i.e., $W(A)$ can be written as a convex function of all minors of $A$, then  the existence of minimizers to \eqref{motivation} was proved by J.M.~Ball  in his pioneering paper \cite{ball77}. We refer, e.g., to \cite{ciarlet,dacorogna} for various results in this direction.  Namely, the existence theory for polyconvex materials can even cope with the important physical assumption, namely,    
\be\label{blowup} W(A)\to+\infty \mbox{ whenever }\det A\to 0_+\ .\ee 
 
On the other hand, there are  \emph{many materials that cannot be modeled by polyconvex stored energies}, prominent examples are materials with microstructure, like shape-memory materials \cite{ball-james1,mueller}. If we give up \eqref{blowup} and suppose that $W$ has  polynomial growth at infinity, e.g. for $c,\tilde c>0$
\be\label{gr-cond}
c(-1+|A|^p)\le W(A)\le \tilde c(1+|A|^p)\ ,
\ee
the existence of a solution to \eqref{motivation} is guaranteed if $W$ is quasiconvex \cite{morrey}, which means that for all $\varphi\in W^{1,\infty}_0(\O;\R^n)$ and all $A\in\R^{n\times n}$   it holds that 
\be\label{quasiconvexity} |\O|W(A)\le\int_\O W(A+\nabla\varphi(x))\,\md x\ .\ee In \eqref{gr-cond} and it the sequel, $|A|^2:=\sum_{i,j=1}^n A_{ij}^2=\sum_{i=1}^n\lambda^2_i(A)$ for $A\in\R^{n\times n}$ with entries $A:=(A_{ij})$ and $\lambda_1(A),\ldots,\lambda_n(A)$ are non decreasingly ordered singular values of $A$, i.e., eigenvalues of $\sqrt{A^\top A}$.

Yet, quasiconvexity is a very complicated property difficult to verify in many cases. Moreover, stored energy densities of materials with microstructure \emph{do not possess this property either}. As a result, solutions to \eqref{motivation} might not exist. Various relaxation techniques were developed \cite{dacorogna,mueller,r} to overcome this drawback for integrands satisfying \eqref{gr-cond}.  Some relaxation results for the case  $W(A)\to +\infty$ if $\det A\to 0$ but $W(A)<+\infty$ even if $\det A<0$ were recently stated in \cite{hm}. In both situations one replaces the integrand by its quasiconvex envelope (the pointwise supremum of all quasiconvex functions not greater than  $W$).     

Another approach  used in variational problems where the integrand satisfies \eqref{gr-cond}  is to  \emph{extend the notion of solutions from Sobolev mappings to parameterized measures}  \emph{called Young measures} \cite{ball3,fonseca-leoni,pedregal,r,tartar,tartar1,y}.  The idea is to describe limit behavior of $\{J(y_k)\}_{k\in\N}$ along a minimizing sequence $\{y_k\}_{k\in\N}$. Nevertheless, the \emph{growth condition \eqref{gr-cond} is still a key ingredient} in these considerations. 

Our goal is to tailor the Young-measure relaxation to functions satisfying \eqref{blowup}. In order to reach this,  we allow $W$ to depend on the inverse of its argument, more precisely, we suppose that $W$ is continuous on invertible matrices and that there exist positive constants $c, \tilde{c} > 0$ such that
 \be\label{ass-1}
 c(-1+|A|^p+|A^{-1}|^p)\le W(A)\le \tilde c(1+|A|^p+|A^{-1}|^p)\ .
 \ee
  
Notice that \eqref{ass-1} implies \eqref{blowup} and that $W$ has  polynomial growth in $|A|$ and $|A^{-1}|$ at infinity.
In the context of nonlinear elasticity, $A$ plays the role of a deformation gradient measuring deformation strain. Then $A^{-1}$ is just another strain measure. We refer, e.g.,~to \cite{c-r, silhavy} for the so-called Seth-Hill family of strain measures or to \cite{haupt} where the physical meaning of the Piola tensor and of the Finger tensor depending on  $A^{-1}A^{-\top}$ and on $A^{-\top}A^{-1}$, respectively, is discussed in great detail. 
Notice that if $y:\Omega\to\R^n$ is a deformation of the reference domain $\O\subset\R^n$  and $y^{-1}:y(\O)\to\O$  is its differentiable inverse then for $x\in\O$
$(\nabla y(x))^{-1}=\nabla y^{-1}(z)$, $z:=y(x)$. Hence, if we exchange the role of the reference and deformed configurations our model 
requires the same integrability for the original deformation gradient as well as for the deformation gradient of the inverse deformation.  
Moreover, if $p\ge 2$ 
\be \sqrt{n}|\O|\le \int_\O|\nabla y(x)||(\nabla y(x))^{-1}|\,\md x\le C\|\nabla y\|_{L^p}\|(\nabla y)^{-1}\|_{L^p}\ ,\ee
which means that if $\|(\nabla y)^{-1}\|_{L^p}\le C $ then $\sqrt{n}|\O|/C\le \|\nabla y\|_{L^p}$ for some $C>0$.  Thus, as a consequence we 
get 
$$
0<\sqrt{n}|\O|/C\le \|\nabla y\|_{L^p}\le C\ .$$

In particular, if $p=+\infty$ and $y\in C^1(\O;\R^n)$ then $y$ is bilipschitz, i.e., $y$ as well as $y^{-1}$ are both Lipschitz maps defining homeomorphism
(called sometimes ``lipeomorphism'') between $\O$ and $y(\O)$. It is a trivial observation, however, that smoothness of $y$ is essential to define lipeomorphism 
and that a positive lower bound on the gradient norm  is not enough to ensure mere invertibility.

This all motivates our idea to perform relaxation in terms of Young measures generated by sequences of matrix-valued mappings $\{Y_k\}_{k\in\N} \subset L^p(\O;\R^{n\times n})$ such that  $\{Y_k^{-1}\}_{k\in\N}\subset L^p(\O;\R^{n\times n})$ is also bounded. We show that, in this case, the Young measures are necessarily supported on invertible matrices and satisfy a certain integral condition, cf. \eqref{def:inverse}. If, additionally, $\det Y_k>0$ almost everywhere in $\O$ for all $k\in\N$ the resulting Young measure is supported on matrices with positive determinant, cf. Theorem~\ref{THM2} and Proposition~\ref{det}. Contrary to the general theory of Young measures generated by $L^p$-maps \cite{schonbek,r}, where 
only the behavior of test functions at infinity is important,  Young measures supported on invertible matrices are also sensitive to the asymptotics  of test functions as the argument approaches a singular matrix.  \emph{However, they allow for a larger class of test functions, namely, those with growth specified in \eqref{ass-1}.} In particular, our test functions are not necessarily continuous on $\R^{n\times n}$.  The precise condition is stated in Theorems~\ref{THM1}, \ref{THM2}. We refer to \cite{kalamajska} for another refinement of Young measures involving discontinuous integrands.  

Although the characterization of Young measures generated by \emph{vector-valued} mappings $\{Y_k\}_{k \in \N} \subset L^p(\O;\R^{n})$ and such that  $\{|Y_k|^{-1}\}_{k\in\N}\subset L^p(\O)$, with $\O$ an interval in $\R$ is already due to Freddi and Paroni \cite{Freddi1}, our manuscript presents, to the authors' knowledge, the first explicit and complete  characterization of Young measures generated by \emph{matrix-valued} mappings $\{Y_k\}_{k \in \N} \subset L^p(\O;\R^{n\times n})$ such that $\{Y_k^{-1}\}_{k\in\N}\subset L^p(\O;\R^{n\times n})$ is also bounded. Moreover, we examine the support of the generated Young measures. 

Also, we completely and and explicitly describe \emph{Young measures generated by gradients} of $W^{1, \infty}(\O, \R^n)$-functions such that the gradients are invertible and the inverse is in $L^{\infty}(\O, \R^n)$, too. The main characterization is exposed in Theorem \ref{THM3}; it is similar to \cite{k-p1}, however, we add additional constraints on the support of the measure as previously in the non-gradient case and introduce a modified  notion of quasi-convexity in \eqref{regular-env}.
The plan of the paper is as follows. After introducing Young measures we state our main results -- Theorems~\ref{THM1}, \ref{THM2}, and \ref{THM3} in Section~\ref{mainresults}. The proofs  of our statements  are left, however, to Section~\ref{mains} after we collect some auxiliary material in Section~\ref{aux}. In particular, Propositions~\ref{Linfty2}, \ref{Linfty1} are of special interest as they form an $L^\infty$ version of our main theorems \ref{THM1}, \ref{THM2}.   We wish to mention that related  results dealing with relaxation for integrands tending to infinity if the determinant approaches zero were proved also in \cite{hm}. Interesting weak* lower semicontinuity theorems for bilipschitz maps are treated in \cite{buliga}. 

Throughout the paper, we use standard notation for Lebesgue $L^p$ and Sobolev  $W^{1,p}$ spaces. We say that $\{u_k\}_{k\in\N}\subset L^1(\O)$ is equi-integrable if we can pick up a subsequence weakly converging in $L^1(\O)$. We refer, e.g.,~to \cite{d-s,fonseca-leoni} for details about equi-integrability  and relative weak compactness in $L^1$. If not said otherwise, $\O\subset\R^n$ is a bounded domain with  a Lipschitz boundary. Finally, $C$ denotes a generic positive constant which may change from line to line.

\subsection{Young measures}
For $p\ge0$ we define the following  subspace of the space
$C(\R^{n\times n})$ of all continuous functions on $\R^{n\times n}$ :
$$
C_p(\R^{n\times n}):=\left\{v\in C(\R^{n\times n}); \lim_{|s|\to\infty}\frac{v(s)}{|s|^p}=0\right\}\ .
$$
Young measures on a bounded  domain $\O\subset\Rn$ are weakly* measurable mappings
$x\mapsto\nu_x:\O\to \rca(\R^{n\times n})$ with values in probability measures;
and the adjective ``weakly* measurable'' means that,
for any $v\in C_0(\R^{n\times n})$, the mapping
$\O\to\R:x\mapsto\A{\nu_x,v}=\int_{\R^{n\times n}} v(s)\nu_x(\d s)$ is
measurable in the usual sense. Let us remind that, by the Riesz theorem,
$\rca(\R^{n\times n})$, normed by the total variation, is a Banach space which is
isometrically isomorphic with $C_0(\R^{n\times n})^*$, where $C_0(\R^{n\times n})$ stands for
the space of all continuous functions $\R^{n\times n}\to\R$ vanishing at infinity.
Let us denote the set of all Young measures by ${\cal Y}(\O;\R^{n\times n})$. It
is known that ${\cal Y}(\O;\R^{n\times n})$ is a convex subset of $L^\infty_{\rm
w}(\O;\rca(\R^{n\times n}))\cong L^1(\O;C_0(\R^{n\times n}))^*$, where the subscript ``w''
indicates the aforementioned property of weak* measurability. Let $S\subset\R^{n\times n}$ be a compact set. A classical result
\cite{tartar,warga} is that for every sequence $\{Y_k\}_{k\in\N}$
bounded in $L^\infty(\O;\R^{n\times n})$ such that $Y_k(x)\in S$ there exists its subsequence (denoted by
the same indices for notational simplicity) and a Young measure
$\nu=\{\nu_x\}_{x\in\O}\in{\cal Y}(\O;\R^{n\times n})$ satisfying
\be\label{jedna2}
\forall v\in C(S):\ \ \ \ \lim_{k\to\infty}v\circ Y_k=v_\nu\ \ \ \
\ \ \mbox{ weakly* in }L^\infty(\O)\ ,
\ee
where $[v\circ Y_k](x)=v(Y_k(x))$ and
\be\label{vnu}
v_\nu(x)=\int_{\R^{n\times n}}v(s)\nu_x(\d s)\ .
\ee
Moreover, $\nu_x$ is supported on $S$  for almost all $x\in\O$.  On the other hand, if $\mu=\{\mu_x\}_{x\in\O}$, $\mu_x$ is supported on $S$ for almost all $x\in\O$ and $x\mapsto\mu_x$ is weakly* measurable then there exist a sequence  $\{Z_k\}_{k\in\N}\subset L^\infty(\O;\R^{n\times n})$, $Z_k(x)\in S$  and \eqref{jedna2} holds with $\mu$ and $Z_k$ instead of $\nu$ and $Y_k$, respectively.

Let us denote by ${\cal Y}^\infty(\O;\R^{n\times n})$ the
set of all Young measures which are created in this way, i.e., by taking
all bounded sequences in $L^\infty(\O;\R^{n\times n})$. Moreover, we denote by ${\cal GY}^\infty(\O;\R^{n\times n})$ the subset of ${\cal Y}^\infty(\O;\R^{n\times n})$ consisting of measures generated by gradients of $\{y_k\}\subset W^{1,\infty}(\O;\R^n)$, i.e., if $Y_k:=\nabla y_k$ in \eqref{jedna2}. Notice that we can suppose that $v\in C(\R^{n\times n})$ in \eqref{jedna2} due to Tietze's theorem.  
A generalization of the $L^\infty$  result was formulated by
Schonbek \cite{schonbek} (cf. also \cite{ball3}): if
$1\le p<+\infty$ then for every sequence
$\{Y_k\}_{k\in\N}$ bounded in $L^p(\O;\R^{n\times n})$ there exists its
subsequence (denoted by the same
indices) and a Young measure
$\nu=\{\nu_x\}_{x\in\O}\in{\cal Y}(\O;\R^{n\times n})$ such that
\be\label{young}
\forall v\in C_p(\R^{n\times n}):\ \ \ \ \lim_{k\to\infty}v\circ Y_k=v_\nu\
\ \ \ \ \ \mbox{ weakly in }L^1(\O)\ .\ee
We say that $\{Y_k\}_{k\in\N}$ generates $\nu$ if (\ref{young}) holds. 
Let us denote by ${\cal Y}^p(\O;\R^{n\times n})$ the set of all Young measures which are obtained through the latter procedure, i.e., by taking all bounded sequences in $L^p(\O;\R^{n\times n})$. It was shown in \cite{kruzik-roubicek} that  if $\nu\in{\cal Y}(\O;\R^{n\times n})$ satisfies the bound 
\be\label{L1cond}
\int_\O\int_{\R^{n\times n}}|s|^p\nu_x(\md s)\,\md x<+\infty  
\ee
then $\nu\in {\cal Y}^p(\O;\R^{n\times n})$.

\bigskip

\section{Main results}
\label{mainresults}
Let $\reg$ denote the set of invertible matrices in $\R^{n\times n}$ and $\regpos$ the set of matrices in  $\R^{n\times n}$ with positive determinant. We write $\inv$ for the continuous  function defined on $\reg$  by $\inv(s):=s^{-1}$, i.e., creating the inversion.
Further, we denote by ${\cal Y}^{p,-p}(\O;\R^{n\times n})$ and by ${\cal Y}^{p,-p}_+(\O;\R^{n\times n})$ the  following subsets of ${\cal Y}^{p}(\O;\R^{n\times n})$:
\be\label{def:inverse}
{\cal Y}^{p,-p}(\O;\R^{n\times n}):=\Big\{\nu\in{\cal Y}^{p}(\O;\R^{n\times n});\ \int_\Omega \int_{\reg} (|s|^p + |s^{-1}|^{p}) \nu_x(\d s)\,\d x<+\infty \ ,\nonumber\\ 
 \nu_x(\reg)=1\mbox{ for a.a. $x\in\O$}\Big \}\ ,
\ee 
 \be\label{def:inverse+}
{\cal Y}^{p,-p}_+(\O;\R^{n\times n}):=\left\{\nu\in{\cal Y}^{p,-p}(\O;\R^{n\times n});\  \nu_x(\regpos)=1\mbox{ for a.a. $x\in\O$}\right\}\ ,
\ee 
and the following subspace of the space of  continuous functions on $\reg$ 
\be
C_{p,-p}(\reg):=\left\{v\in C(\reg); \lim_{|s|+|s^{-1}|\to\infty}\frac{v(s)}{|s|^p+|s^{-1}|^p}=0\right\}\ .
\ee
Our main results are summarized in the following theorems.

\bigskip

\begin{theorem}\label{THM1}
Let $+\infty>p\ge 1$, let $\O\subset\R^n$ be open and bounded, and let $\{Y_k\}_{k\in\N}$, $\{Y^{-1}_k\}_{k\in\N} \subset L^p(\O;\R^{n\times n})$ be bounded.  Then there is a subsequence of  $\{Y_k\}_{k\in\N}$ (not relabeled)  and $\nu\in {\cal Y}^{p,-p}(\O;\R^{n\times n})$ such that for every $g\in L^\infty(\O)$ and every $v\in C_{p,-p}(\reg)$ it holds that
\be\label{mainpassage}
\lim_{k\to\infty} \int_\O v(Y_k(x))g(x)\,\md x=\int_\O\int_\reg v(s)\nu_x(\md s)g(x)\,\md x\ ,
\ee
Conversely, if $\nu\in{\cal Y}^{p,-p}(\O;\R^{n\times n})$ then there is a bounded sequence $\{Y_k\}_{k\in\N}\subset  L^p(\O;\R^{n\times n})$ such that $\{Y_k^{-1}\}_{k\in\N}\subset  L^p(\O;\R^{n\times n})$ is also bounded and \eqref{mainpassage} holds for all  $g$ and $v$ defined above.
\end{theorem}

\bigskip

\begin{theorem}\label{THM2}
Let $+\infty>p\ge 1$, let $\O\subset\R^n$ be open and bounded, and let $\{Y_k\}_{k\in\N}$, $\{Y^{-1}_k\}_{k\in\N} \subset L^p(\O;\R^{n\times n})$ be bounded and for every $k\in\N$ $\det Y_k>0$ almost everywhere in $\O$.  Then there is a subsequence of  $\{Y_k\}_{k\in\N}$ (not relabeled)  and $\nu\in {\cal Y}^{p,-p}_+(\O;\R^{n\times n})$ such that for every $g\in L^\infty(\O)$ and every $v\in C_{p,-p}(\reg)$ \eqref{mainpassage} holds.

Conversely, if $\nu\in{\cal Y}^{p,-p}_+(\O;\R^{n\times n})$ then there is a bounded sequence $\{Y_k\}_{k\in\N}\subset  L^p(\O;\R^{n\times n})$ such that $\{Y_k^{-1}\}_{k\in\N}\subset  L^p(\O;\R^{n\times n})$ is also bounded, for every $k\in\N$ $\det Y_k>0$ almost everywhere in $\O$,  and \eqref{mainpassage} holds for all  $g$ and $v$ defined above.
\end{theorem}

\bigskip

Notice that  $C(\R^{n\times n})\cap C_{p,-p}(\reg)=C_p(\R^{n\times n})$, so we allow for a larger class of test function in Theorems~\ref{THM1},~\ref{THM2} compared with
the original result by Schonbek \cite{schonbek} mentioned in \eqref{young}. In particular, our test functions are not necessarily continuous on the whole $\R^{n\times n}$.

\begin{remark}
For simplicity, we formulated Theorems \ref{THM1} and \ref{THM2} as well as Definitions \eqref{def:inverse}, \eqref{def:inverse+} \emph{symmetrically} in $p$ in the sense that both the generating sequence as well as its inverse are bounded in $L^p(\Omega; \R^{n \times n})$. We could, however, also define for $\infty>p,q\ge 1$ 
\be
{\cal Y}^{p,-q}(\O;\R^{n\times n}):=\Big\{\nu\in{\cal Y}^{p}(\O;\R^{n\times n});\ \int_\Omega \int_{\reg} (|s|^p + |s^{-1}|^{q}) \nu_x(\d s)\,\d x<+\infty\ ,\nonumber\\  
\nu_x(\reg)=1\mbox{ for a.a. $x\in\O$}\Big\}\ ,
\ee 
 \be
{\cal Y}^{p,-q}_+(\O;\R^{n\times n}):=\left\{\nu\in{\cal Y}^{p,-q}(\O;\R^{n\times n});\  \nu_x(\regpos)=1\mbox{ for a.a. $x\in\O$}\right\}\ .
\ee
Then Theorems \ref{THM1} and \ref{THM2} hold with the single modification that $\left\{Y_k^{-1}\right\}_{k\in\N}$ is bounded $L^q(\Omega; \R^{n \times n})$.
\end{remark}

\begin{remark}
We could also define the set
\be
{\cal Y}^{p,f(\cdot)}(\O;\R^{n\times n}):=\Big\{\nu\in{\cal Y}^{p}(\O;\R^{n\times n});\ \int_\Omega \int_{\reg} (|s|^p + f(s^{-1})) \nu_x(\d s)\,\d x<+\infty \ ,\nonumber\\ 
 \nu_x(\reg)=1\mbox{ for a.a. $x\in\O$}\Big \}\ ,
\ee 
with $f(\cdot) \geq |\det(\cdot)|^q$ for some $q > 0$. Obvious modifications of the proofs below give that $\nu$ is in ${\cal Y}^{p,f(\cdot)}(\O;\R^{n\times n})$ if and only if it can be generated by a sequence of invertible matrices inverses $\left\{Y_k^{-1}\right\}_{k\in\N}$ bounded in  $L^q(\Omega; \R^{n \times n})$. Defining this set allows us us relaxe even a larges class of functions than  $C_{p,-p}(\reg)$. 
\end{remark}

\bigskip

The next result shows that the weak limit  of a sequence of gradients with positive determinants inherits this property if we control the behavior of the inverse.   

\bigskip

\begin{proposition}\label{det}
Let $p>n$. If $y_k\wto y$ in $W^{1,p}(\O;\R^n)$, for all $k\in\N$ $\det\nabla y_k>0$ a.e.~in $\O$, and $\{(\nabla y_k)^{-1}\}\subset L^p(\O;\R^{n\times n})$ is bounded   then $\det\nabla y>0$ a.e.~in $\O$. Moreover, every Young measure generated by a subsequence of $\{\nabla y_k\}_{k\in\N}$ is supported on $\regpos$. 
\end{proposition}

\bigskip

We now turn to a characterization of  gradient Young measures supported on invertible matrices. This allows us to formulate new relaxation and weak* lower semicontinuity theorems for integrands continuous on invertible matrices (or at least a compact subset of them) but not necessarily well-defined at singular matrices and which tend to infinity if their argument approaches a non-invertible matrix, provided their infimizing sequence is uniformly bounded in an appropriate sense. We shall see that the characterization is analogous to the one obtained by Kinderlehrer and Pedregal for gradient Young measures \cite{k-p1,k-p}, however, the quasiconvex envelope is replaced by possibly  a more restrictive one and a condition on the support of the Young measure is added.

We will define the following sets of Young measures generated by bounded and invertible  gradients of $W^{1,\infty}(\O;\R^n)$ maps:
\be
\mathcal{GY}^{+\infty,-\infty}_\varrho(\O;\R^{n\times n}):=\Big\{\nu\in{\cal Y}^\infty(\O;\R^{n\times n});\ \exists \{y_k\}\subset W^{1,\infty}(\O;\R^n)\ ,\ \\ \nonumber \mbox{for a.a.~$x\in\O$ } \{\nabla y_k(x)\}\subset\Rrho\mbox{ and $\{\nabla y_k\}$ generates $\nu$}\Big\}\ 
\ee 
with $\Rrho:=\{A\in\reg;\ \max(|A|,|A^{-1}|)\le\varrho\}$ and $\Grad:=\cup_{\varrho>0}\mathcal{GY}^{+\infty,-\infty}_\varrho(\O;\R^{n\times n})$.

Before stating our characterization of gradient Young measures generated by invertible gradients we will need the following definitions. 
\bigskip

Put $R^{n\times n}_{+\infty}:=\reg$ and denote  for $\varrho\in(0;+\infty]$   

$$\mathcal{O}(\varrho):=\{ v:\R^{n\times n}\to\R\cup\{+\infty\};\ v\in C(\Rrho)\ ,\ v(s)=+\infty \mbox { if $s\in\R^{n\times n}\setminus\Rrho$}\}\ .\ $$

If $F\in\R^{n\times n}$ and $v\in\mathcal{O}(\varrho) $ we denote by $Q_{\rm inv}v:\R^{n\times n}\to\R\cup\{+\infty\}$  the 
function 
\be\label{regular-env}
Q_{\rm inv}v(F):=|\O|^{-1}\inf_{y\in\mathcal{U}_F}\int_\O v(\nabla y(x))\,\md x\ ,
\ee 
where 
\be\mathcal{U}_F:=\{y\in W^{1,\infty}(\O;\R^n); (\nabla y)^{-1}\in L^\infty(\O;\R^{n\times n}),\ \ y(x)=Fx \  {\rm for  }\ x\in\partial \Omega\}\ . \label{UF}\ee 
\bigskip

\begin{remark}
Note that the boundeness of the inverse in $L^\infty(\O;\R^{n\times n})$ of the test functions demanded in the defintion of $Q_{\rm inv}$ in \eqref{regular-env} actually means that $|\det(\nabla v)|\geq c > 0$. In fact therefore, one could consider instead of \eqref{qc0} in Theorem \ref{THM3} the standard Jensen inequality as in \cite{k-p1} with this constraint included.
\end{remark}
\begin{remark}

It is not obvious that $\mathcal{U}_F$ is non-empty if $F$ is singular. If $\mathcal{U}_F$ were empty $Q_{\rm inv}v(F)$ would be equal to $+\infty$ even if $v\in C(\reg)$  and hence $Q_{\rm inv}v(F)$ would not be always finite. However, in Remark \ref{finiteness} it is noted that this situation will not occur.

Similarly as for quasiconvexity, one can show that $Q_{\rm inv}v$ does not depend on the Lipschitz domain $\O$ used in its definition.
\end{remark}

\bigskip

\begin{theorem}\label{THM3}
Let $\O\subset\R^n$ be an open bounded Lipschitz domain and let $\nu\in{\cal Y}^\infty(\O;\R^{n\times n})$. Then $\nu\in \Grad$ if and only if the following three conditions hold
 \be\label{supp}
{\rm supp }\,\nu_x\subset\Rrho \mbox{ for a.a. $x\in\O$ and some $\varrho>0$}, 
\ee 
 \be\label{firstmoment0}
 \exists\ u\in W^{1,\infty}(\O;\R^n)\ :\  \nabla u(x)=\int_{\reg} s \nu_x(\d s)\  \ ,\ee
for a.a. $x\in\O$ all $\tilde\varrho\in(\varrho;+\infty]$, and all  $v\in \mathcal{O}(\tilde\varrho)$  the following  inequality
is valid \be\label{qc0} Q_{\rm inv}v(\nabla u(x))\le\int_{\reg} v(s)\nu_x(\md s)\ .\ee

\end{theorem}

\bigskip

\begin{remark}
Clearly, \eqref{qc0} holds for all $v\in C(\R^{n\times n})$,  too. Indeed,  choose $\tilde\varrho$ as in the theorem and define 
$$
\tilde v(s):=\begin{cases}
v(s) &\mbox{ if $s\in R^{n\times n}_{\tilde\varrho}$,}\\
+\infty &\mbox{ otherwise.}
\end{cases}
$$
Then $\tilde v\ge v$ and $\tilde v\in \mathcal{O}(\tilde\varrho)$. Hence, we have from \eqref{qc0} and  \eqref{regular-env}
$$\int_{\reg} v(s)\nu_x(\md s)=\int_{\reg} \tilde v(s)\nu_x(\md s)\ge Q_{\rm inv}\tilde v(\nabla u(x))\ge Q_{\rm inv}v(\nabla u(x))\ .$$

\end{remark}

\bigskip

\section{Auxiliary results}
\label{aux}

Let us start by recalling the definition
\begin{align}
\label{Rvarrho} \Rrho&:=\{A\in\reg;\ \max(|A|,|A^{-1}|)\le\varrho\}\ , \\
\intertext{and defining analogously}
\label{Rvarrho+ }\!\!\!\! \RrhoPlus&:=\{A\in \Rrho;\  {\rm det}\ A>0\}\ .
\end{align}
Then the following holds:
\begin{lemma}
$\Rrho$ is compact in $\R^{n\times n}$ for every $\varrho>0$. Moreover, the set $\RrhoPlus$ is also compact for every $\varrho>0$.
\end{lemma}
\bigskip
\noindent{\it Proof.}
Clearly, $\Rrho$ is bounded. Consider a sequence $\{A_k\}_{k\in\N}\subset \Rrho$ such that $A_k\to A$. We must show that $A\in \Rrho$.  If $\det A=0$ then by continuity  $\det A_k\to 0$ as $k\to\infty$ and due to the bound
$|\det B|\le C|B|^n$, $C>0$ for all $B\in\Rm$ we would have 
$$|1/\det A_k|=|\det A_k^{-1}|\le C|A_k^{-1}|^n\to\infty\ .$$
Hence, $A_k\not\in \Rrho$ if $k\ge k_0$ which is a contradiction. Therefore, $A\in\reg$. The continuity of the matrix inverse $A_k^{-1}\to A^{-1}$ and $|A_k^{-1}|\to |A^{-1}|$ yields, in consequence, that $A\in \Rrho$ and $\Rrho$ is bounded and closed. Compactness of $\RrhoPlus$ follows by  continuity of the function $A\mapsto \det A$.  \hfill $\Box$

\bigskip

\begin{remark}
We have $\reg=\cup_{\varrho\in\N}\Rrho$ and $\regpos=\cup_{\varrho\in\N}\RrhoPlus$, i.e., the open sets $\reg$ and $\regpos$ are both  $\sigma$-compact. 
\end{remark}

\bigskip
\noindent For every $v:\R^{n\times n}\to\R$  we define $\hat v:\R^{n\times n}\to\R$:
\be\label{hatv}
\hat v(s):=
\begin{cases}
v(s^{-1}) & \mbox{ if $s\in\reg$, }\\
0 & \mbox{ otherwise.}
\end{cases}
\ee

\bigskip

\noindent We define the following subspace of $C_0(\R^{n\times n})$:
\be\label{c0reg}
\C:=\{v\in C_0(\R^{n\times n});\ v(s)=0\mbox { if } \det s=0\}\, \ee 
equipped with the supremum norm. 
Notice, that  $v=\hat{\hat v}$ for every $v\in\C$.

\bigskip

\begin{lemma}
\label{PropOfC}
$\C$ is a \emph{separable Banach space} with respect to the standard maximum norm for continuous functions. Moreover, $\C = \overline{\bigcup_{\varrho > 0} C_{\Rrho}(\Rm)}$, where $C_{\Rrho}(\R^{n \times n}): = \{ \varphi \in C_0(\R^{n \times n}), \, \mathrm{supp}\, \varphi \subset \Rrho\}$. 
\end{lemma}

\begin{proof}
First of all, notice that $\C$ is closed in $C_0(\R^{n \times n})$. Indeed, take a sequence $\{\phi_k\}_{k\in\N} \subset \C$ such that $\phi_k \to \phi$ in $C_0(\R^{n \times n})$. Then, in particular, $\phi_k(A) \to \phi(A)$ for all $A \in \R^{n \times n}$. Hence, also $\phi(A) = 0$ for every singular $A$ meaning that $\phi \in \C$. Therefore, $\C$ is also a Banach space.

Clearly, any $\phi$ in $C_{\Rrho}(\R^{n \times n})$ is also in $\C$ for any $\varrho > 0$. Hence also $ {\bigcup_{\varrho > 0} C_{\Rrho}(\R^{n \times n})} \subset \C$ and, because $\C$ is closed, the same holds for the closure. On the other hand, take $\phi\in \C$ and define for every $\varrho$ a \emph{smooth cut-off function $\Phi_\varrho$} which is 1 on $\Rrho$ and 0 on $ \reg\setminus R^{n\times n}_{\varrho+1}$. (Note that $\Phi_\varrho$ can be found as follows: Define $\Theta_\varrho$ a smooth function which is 1 inside the ball $B(0,\varrho) \subset \Rm$ and equals 0 on $ \Rm\setminus B(0,\varrho+1)$. Now we may set $\Phi_\varrho(s): = \Theta_\varrho(s) \hat\Theta_\varrho (s)$. Note that, since $B(0,\varrho+1)$ is a strict subset of $\Rm$, $\Phi_\varrho$  is indeed smooth.) Then $\phi$ can be approximated by the set of functions $\{\phi \cdot \Phi_\varrho\}_{\varrho>0}$ if we can show that 
\begin{equation}
\forall \epsilon > 0 \, \, \exists \varrho_0 >0: \quad |\phi(A)| < \epsilon \, \quad \forall A \in \R^{n \times n} \setminus R^{n\times n}_{\varrho}.
\label{approxC0}
\end{equation}
To see this, suppose for contradiction that \eqref{approxC0} does not hold and that there exists $\epsilon > 0$ and  $\{A_{\varrho}\} \subset \R^{n\times n} \setminus \Rrho$ such that $\phi(A_{\varrho}) \ge \epsilon$. Clearly, $\{A_{\varrho}\}$ must be bounded since $\phi \in C_0(\Rm)$. Therefore, pick a subsequence of $ \{A_\varrho\}$ (not relabeled) such that $A_\varrho \to A$. Then, also $\phi(A) \geq \epsilon$ from which it follows that $A$ is  an invertible matrix. But since $\bigcup_{\varrho \geq 0} \Rrho = \reg$ there has to exist  $\tilde{\varrho}$ such that $A \in R^{n\times n}_{\tilde{\varrho}}$. Yet, from the construction, for $\varrho$ large enough $A_\varrho$ are \emph{not} elements of $R_{\tilde{\varrho}}^{n \times n}$, a contradiction. 

For the separability we use the classical result that subspaces of separable metric spaces are again separable \cite{lohman}.
\end{proof}

Notice that if $v\in\C$ then  $\hat v\in\C$.  Indeed, if $s_0\in\reg$ then there is a $\delta$-neighborhood of $s_0^{-1}$, $B(s_0^{-1},\delta)$,  such that $B(s^{-1}_0,\delta)\subset\reg$. The function $v$ is continuous on  $B(s^{-1}_0,\delta)$, so for every $\varepsilon>0$ we have $|\hat v(s_0)-\hat v(s)|=|v(s_0^{-1})-v(s^{-1})|\le\varepsilon$ if $\delta>0$ is small enough.   If $s_0$ is singular, then $\hat v(s_0)=0$ and  $|s^{-1}|$ is arbitrarily large on $B(s_0,\delta)\cap\reg$. Hence, $|\hat v(s)|=|v(s^{-1})|<\varepsilon$ on $B(s_0,\delta)\cap\reg$ if $\delta$ is small. On the other hand, for singular matrices $s_0 \in  B(s_0,\delta)$ $\hat v(s_0)=0$, anyway.

The following lemma is a simple but useful observation. Namely, it asserts that it is enough to test by functions from the separable space $\C$  to identify equal measures supported on $\reg$.

\bigskip

\begin{lemma}\label{equalityofmeasures} Let $\nu,\mu\in\rca(\R^{n\times n})$ and let both be supported on $\reg$. If for every $v\in\C$
\be\label{eqonreg}\int_{\reg}v(s)\nu(\md s)=\int_{\reg}v(s)\mu(\md s)\ ,
\ee
then $\nu=\mu$, i.e., \eqref{eqonreg} holds even for all $v\in C_0(\R^{n\times n})$.

\end{lemma}

\noindent
{\it Proof.} Take $v\in C_0(\R^{n\times n})$.
Define, similarly as in the proof of Lemma \ref{PropOfC}, for every $\varrho$ the \emph{smooth cut-off function $\Phi_\varrho$} which is 1 on $\Rrho$ and 0 on $\reg\setminus R^{n\times n}_{\varrho+1}$ and $v_\varrho(s):=v(s)\Phi_\varrho(s)$ for all $s\in\R^{n\times n}$. Then $v_\varrho \in \C$ and $|v_\varrho|\le |v|$. The proof is finished by the Lebesgue dominated convergence theorem for which we notice that $v_\varrho \to v$ pointwise everywhere on $\reg$.
\hfill
$\Box$
\bigskip

\begin{proposition}\label{Linfty2}
Let $\nu\in\mathcal{Y}(\O;\R^{n\times n})$ and suppose that there is $\varrho>0$ such that  for almost all $x\in\O$ supp\ $\nu_x\subset \Rrho$. 
Then  there exists $\{Y_k\}_{k\in\N}\subset L^\infty(\O;\R^{n\times n})$  such that $\{Y_k(x)\}_{k\in\N}\subset \Rrho$ for almost all $x\in\O$ and $\{Y_k\}_{k\in\N}$ generates $\nu$.
\end{proposition}
\bigskip
\noindent{\it Proof.} This is a classical result mentioned in \eqref{jedna2}. See e.g.~\cite[Th.~1]{tartar} for details.
\hfill
$\Box$

\bigskip

\begin{remark}\label{plus}
Proposition~\ref{Linfty2} still holds if we replace $\Rrho$ by $\RrhoPlus$ which is compact, as well.
\end{remark}

\bigskip

\begin{proposition}\label{Linfty1}
Let $\varrho>0$ and let $\{Y_k\}\subset L^\infty(\O;\R^{n\times n})$, $\{Y_k\}\subset \Rrho$  for almost all $x\in\O$ and all $k\in\N$. 
If $\{Y_k\}$ generates $\nu\in{\cal Y}^\infty(\O;\R^{n\times n})$ and if $\{Y_k^{-1}\}$ generates $\mu\in{\cal Y}^\infty(\O;\R^{n\times n})$ then for almost all $x\in\O$ and every continuous $f:\Rrho\to\R$  it holds
\be\label{munurelation}
\int_{\Rrho}f(s)\mu_x(\md s)=\int_{\Rrho} \hat{f}(s)\nu_x(\md s)\ .
\ee
Moreover, $\mathrm{supp} \, \nu_x \subset \Rrho$ for almost all $x \in \Omega$.
\end{proposition}
\bigskip
\noindent{\it Proof.}
First of all, recall that \cite{balder,valadier} for almost all $x\in\O$ $\nu_x$ is supported on the set $\cap_{l=1}^\infty\overline{\{Y_k(x);\ k\ge l\}}$, i.e., $\nu_x$ is supported on $\Rrho$. 
Further, notice that $\{Y^{-1}_k(x)\}\subset \Rrho$ for a.a. $x \in \Omega$. If $f:\Rrho\to\R$ is continuous, so is $F:\Rrho\to\R$, $F(s):=f(s^{-1})$. Then we have for any $g\in L^1(\O)$
$$\lim_{k\to\infty} \int_\O f(Y^{-1}_k(x))g(x)\,\md x=\int_\O\int_{\Rrho}f(s)\mu_x(\md s)g(x)\,\md x\ .$$
At the same time,
$$\lim_{k\to\infty} \int_\O F(Y_k(x))g(x)\,\md x=\int_\O\int_{\Rrho}F(s)\nu_x(\md s)g(x)\,\md x=\int_\O\int_{\Rrho}f(s^{-1})\nu_x(\md s)g(x)\,\md x\ .$$
\hfill $\Box$
\bigskip
\begin{proposition}\label{support}
Let $\{Y_k\}_{k\in\N}\subset L^p(\O;\R^{n\times n})$ generate $\nu\in\mathcal{Y}^p(\O;\R^{n\times n})$  and let $\int_\Omega |{\rm det} Y_k^{-1}|^q \d x \leq C$ for some $C > 0$ and some $q >0$.
Then for almost all $x\in\O$ $\nu_x$ is supported on $\reg$ in the sense that $\nu_x(\R^{n\times n}\setminus\reg)=0$ for almost all $x\in\O$.
Moreover, if even $\mathrm{det}\, Y_k > 0$ a.e. in $\O$ then $\nu_x$ is for almost all $x \in \O$ supported (in the above sense) on the set of invertible matrices with positive determinant.

\end{proposition}
\bigskip
\noindent{\it Proof.}
 Assume that the first assertion did not hold, i.e., that there existed a measurable $\omega \subset \O$ with positive measure such that $\int_\omega \int_{\R^{n\times n}\setminus\reg} \nu_x (\d s) \d x > 0$. Then for any $\epsilon >0$ define a smooth cut-off $\Phi^{\mathrm{det},0}_\epsilon$ such that $\Phi^{\mathrm{det},0}_\epsilon(s) = 1$ on $\R^{n\times n}\setminus\reg$ and $\Phi^{\mathrm{det},0}_\epsilon(s) = 0$ for all $s \in \reg$ such that $|\det s| \geq \epsilon$ ($\Phi^{\mathrm{det},0}_\epsilon$ can be found as follows: first of all find a smooth $\varphi_\epsilon: \R \to \R$ such that $\varphi_\epsilon (0) = 1$ and $\varphi_\epsilon(x) = 0$ for $|x| > \epsilon$. Then define $\Phi^{\mathrm{det},0}_\epsilon(s) = \varphi_\epsilon(\det s))$. We have by \eqref{young}
\be\label{firstpassage}
\lim_{k\to\infty}\int_\O \Phi^{\mathrm{det},0}_\epsilon(Y_k(x))\,\md x=\int_\O\int_{\R^{n\times n}}\Phi^{\mathrm{det},0}_\epsilon(s)\nu_x(\md s)\,\md x\ge\int_\omega \int_{\R^{n\times n}\setminus\reg} \nu_x (\d s) \d x=:\delta>0\ .
\ee
Hence, there is $k_0\in\N$ such that  $\int_\O \Phi^{\mathrm{det},0}_\epsilon(Y_k(x))\,\md x>\delta/2$ if $k>k_0$. This means that there is always a measurable set $\omega(k)\subset\O$,  $|\omega(k)|>\delta/2$   such that 
$|{\rm det}\ Y_k(x)|^q<\varepsilon^q$ if $x\in\omega(k)$. Consequently,  $|{\rm det}\ Y^{-1}_k(x)|^q>\varepsilon^{-q}$ if $x\in\omega(k)$. 
Thus, for every $k>k_0$ 
\be\label{secondpassage}\int_\O|{\rm det}\ Y_k^{-1}(x)|^q\,\md x\ge \int_{\omega(k)}|{\rm det}\ Y_k^{-1}(x)|^q\,\md x\ge \frac{\delta}{2\varepsilon^q}\ .
\ee
As $\varepsilon>0$ is arbitrary it contradicts the bound $\int_\Omega |{\rm det} Y_k^{-1}|^q \d x \leq C$.

As to the second assertion we proceed analogously only we define instead of $\Phi^{\mathrm{det},0}_\epsilon$ the smooth cut-off $\Phi^{\mathrm{det},+}_\epsilon$ which is 1 on all matrices $s$ for which $\det s \leq 0$ and 0 on matrices for which $\det s \geq \epsilon$. Then, if $\nu_x$ was not, for almost all $x \in \O$,  supported on the set of invertible matrices with positive determinant, again there would be a measurable subset of $\O$ with positive measure, such that $\int_\omega \int_{\reg\setminus\regpos} \nu_x (\d s) \d x > 0$ which analogously to \eqref{firstpassage} means that in some set $\omega(k)\subset\O$  $0 \leq \det Y_k \leq \epsilon$.  This yields a contradiction for $\epsilon \to 0$  because of  \eqref{secondpassage}.
\hfill
$\Box$
 
\bigskip

\noindent For further notation, let us denote the set $C^{p,-p}(\reg)$ as
$$
C^{p,-p}(\reg) = \{f\in C(\reg);\  |f(s)|\leq C(1+|s|^p + |s^{-1}|^p) \, \, \forall s \in \reg  \}\ .
$$
\begin{lemma}\label{equality}
Let $\nu \in{\cal Y}^{p,-p}(\O;\R^{n\times n}), \, \mu \in {\cal Y}^{p,-p}(\O;\R^{n\times n})$.
Let $f \in C^{p,-p}(\reg)$ and $\hat f$ be defined as in \eqref{hatv} with $f$ instead of $v$. Let also,
\begin{equation}\label{restriction}
\int_\Omega \int_{\reg} f^\varrho(s) \mu_x(\d s)\, \d x = \int_\Omega \int_{\reg} \hat{f}^\varrho(s) \nu_x (\d s) \d x
\end{equation}
for all $f^\varrho \in C_\Rrho(\Rm)$, for any $\varrho > 0$. Then
\begin{equation}
\int_\Omega \int_{\reg} f(s) \mu_x(\d s)  \d x = \int_\Omega \int_{\reg} \hat f(s) \nu_x (\d s) \d x\ .
\label{InverseMeasureRelationship}
\end{equation}
\end{lemma}

\bigskip

\begin{proof} 
Take any $f \in C^{p, -p}(\reg)$ and define (the same was as in the proof of Lemma \ref{PropOfC} the smooth cut-off $\Phi_\varrho$. Then $f(s) \Phi_\varrho (s) \to f(s)$ pointwise for all $s \in \reg$ and  similarly (from continuity) $\hat f(s) \Phi_\varrho (s) \to \hat f(s)$ as $\varrho\to\infty$. Note also that    
\begin{align*}
\int_\Omega \int_{\reg} |f(s)| \Phi_\varrho(s) \mu_x(\d s) \d x &\leq \int_\Omega C(1+ |s|^p +|s^{-1}|^{p}) \mu_x(\d s) \d x \leq \tilde C \\
\int_\Omega \int_{\reg} |\hat f(s)| \Phi_\varrho (s)\nu_x (\d s) \d x 
&\leq \int_\Omega C(1 +|s^{-1}|^{p} + |s|^p) \nu_x(\d s) \d x \leq \tilde C
\end{align*}
are bounded independently of $\varrho$. Here $\tilde C>0$ is a constant.  Hence we may exploit Lebesgue's dominated convergence theorem  to prove the assertion.  
\end{proof}
 
\bigskip

\begin{proposition}\label{integrability}
Let $p \in [1, \infty)$ and $\{Y_k\}\subset L^p(\O;\R^{n\times n})$, $\{Y_k^{-1}\}\subset L^p(\O;\R^{n\times n})$ be bounded and $\{Y_k(x)\}\subset\reg$ for almost all $x\in\O$. Then there is a subsequence of $\{Y_k\}$ (not relabeled) such that this subsequence generates a Young measure $\nu\in{\cal Y}^{p, -p}(\O;\R^{n\times n})$.

Moreover, if we denoted $\mu$ the Young measure generated by (a further subsequence of) $\{Y_k^{-1}\}$ then \eqref{InverseMeasureRelationship} holds for all $f \in C^{p, -p}(\reg)$. 

\end{proposition}
\bigskip
\noindent{\it Proof.} It follows from \eqref{young} that  a (not relabeled) subsequence of $\{Y_k\}$ generates a Young measure $\nu\in{\cal Y}^p(\O;\R^{n\times n})$ and $\{Y^{-1}_k\}$ generates a Young measure $\mu\in{\cal Y}^p(\O;\R^{n\times n})$. As $\int_\Omega |\det (Y_k^{-1})|^{p/n} \d x \leq  C \int_\Omega |Y_k^{-1}|^{p} \d x \leq C$, we know from Proposition~\ref{support} that $\nu_x$ and $\mu_x$ are both supported on $\reg$ for almost all $x\in\O$.  
We have for  all $g\in L^\infty(\O)$ and all $v\in\C$ 
$$
\lim_{k\to\infty}\int_\O \hat v(Y_k(x))g(x)\,\md x =\int_\O\int_{\reg} \hat v(s)\nu_x(\md s)g(x)\,\md x\ , $$
$$
\lim_{k\to\infty}\int_\O v(Y^{-1}_k(x))g(x)\,\md x =\int_\O\int_{\reg} v(s)\mu_x(\md s)g(x)\,\md x\ , $$
where $\hat{v}$ given by \eqref{hatv} is again in $\C$. This means that for  all $g\in L^\infty(\O)$ and all $v\in\C$
\be
\int_\O\int_{\reg} \hat v(s)\nu_x(\md s)g(x)\,\md x=\int_\O\int_{\reg} v(s)\mu_x(\md s)g(x)\,\md x\ .
\ee 
If we define
$$\int_\reg v(s)\hat\nu_x(\md s):=\int_\reg \hat{v}(s)\nu_x(\md s)\ ,$$
we get by Lemma~\ref{equalityofmeasures} that $\hat\nu=\mu$.

Therefore it remains only to prove that $\int_\Omega \int_\reg (|s|^p+|s^{-1}|^p) \nu_x(\d s) \d x$ is bounded. Boundedness of the first part is guaranteed due to the fact that $\nu\in{\cal Y}^p(\O;\R^{n\times n})$. To see the second part note that $|\cdot|^p \circ\inv$ is a \emph{continuous, bounded from below} in $\reg$ and hence \cite{pedregal}
\be
\int_\O \int_\reg |s^{-1}|^p \nu_x (\d s) \d x &=& \int_\O \int_\reg (|\cdot|^p \circ\inv)(s) \nu_x (\d s) \d x \leq \liminf_{k \to \infty} \int_\O (|\cdot|^p \circ\inv)(Y_k) \d x \nonumber\\
&=& \liminf_{k \to \infty} \int_\O |Y_k^{-1}|^p \d x < C.
\ee

Therefore, by Lemma \ref{equality}, \eqref{InverseMeasureRelationship} holds for all $f \in C^{p, -p}(\reg)$ 

\bigskip

\begin{proposition}\label{explicit}
Let $\nu\in{\cal Y}^{p, -p}(\O;\R^{n\times n})$. Then there is a generating sequence $\{Y_k\}\subset L^p(\O;\R^{n\times n})$ such that $\{Y^{-1}_k\}\subset L^p(\O;\R^{n\times n})$ is bounded. Moreover, $\{|Y^{-1}_k|^p\}$ as well as $\{|Y_k|^p\}$ are equi-integrable. 

\end{proposition}
\bigskip
\noindent{\it Proof.}
Notice, that inevitably for a.a. $x\in\O$ supp $\nu_x\subset\reg$ (cf \eqref{def:inverse}). Therefore, define smooth cut-off functions $\Phi_\varrho$ as in the proof of Lemma \ref{equality} and set
\be\label{newmeasure}
\nu_{x}^{\r}\ =\ \Phi_\varrho\nu_x+\left( \int_{\reg}\left(1-\Phi_\varrho(s)\right)\nu_x(\md s) \right)\delta_I\ ,
\ee 
where $\delta_I$ denotes the Dirac measure supported at the identity matrix. 
It is only a simple observation that, for all $\varrho \in\N$ and a.a.
$x\in\O$, $\nu^{\r}_{x}$ is a probability Radon measure supported on
$R^{n \times n}_{\r+1}$ and that the mapping $x\mapsto \nu^{\r}:\O\to \mathrm{rca}(\Rm)$ is weakly measurable. 
Altogether, we see that $\nu^\r$ defined by \eqref{newmeasure} is a Young
measure, i.e. $\nu^\r\in{\cal Y}(\O;\Rm)$.  We have from Propositions~\ref{Linfty2},\ref{Linfty1} that there is $\{Y_k^\r(x)\}\subset R^{n \times n}_{\r+1}$, $\{(Y_k^\r)^{-1}(x)\}\subset R^{n \times n}_{\r+1}$ for a.a. $x \in \Omega$ such that they generate $\nu^\r$ and $\mu^\r$, respectively, where 
\be
\int_\reg \hat{v}(s) \nu^\varrho_x (\d s) = \int_\reg v(s) \mu^\varrho_x (\d s)
\label{RelMuRNuR}
\ee
holds for all $v \in \C$. 

Now we want to show that, for any $v\in \C$, it holds
$\lim_{\r\to\infty}v_{\nu^\r}=v_\nu$ weakly in $L^1(\O)$, where $v_{\nu_\r}$ is defined again by \eqref{vnu} with $\nu_\r$ in place of $\nu$. Indeed, for any
$g \in L^{\infty}(\O)$ we can estimate
$$
\lim_{\r \to \infty} \int_\Omega g(x) v_{\nu^\r}(x) \d x = \lim_{\r \to \infty} \int_\Omega g(x) \int_\reg v(s) \Phi_\varrho(s) \nu_x (\d s) \d x + \lim_{\r \to \infty} v(I)\int_\Omega g(x) \int_\reg(1-\Phi_\varrho(s)) \nu_x (\d s) \d x
$$
Now, thanks to the proof of Lemma \ref{PropOfC}, we know that $v \Phi_\varrho$ converges strongly in the $C_0$-norm to $v$ and hence the first limit converges to $v_\nu$. As for the second limit $\Phi_\varrho(s)$ converges pointwise to 1 for all $s$ in $\reg$ and therefore, thanks to Lebesgue's dominated convergence theorem (note that $1-\Phi_\varrho$ is bounded by 1 which is $\nu_x- integrable$), $ \int_\reg(1-\Phi_\varrho(s)) \nu_x (\d s)$ converges to 0 for a.a. $x \in \Omega$. Exploiting the dominated convergence once again, the second limit approaches zero. Hence we are in the situation that
$$
\lim_{\r\rightarrow\infty}\lim_{k\rightarrow\infty}v(Y^\r_k)
=v_{\nu}\mbox{ weakly in }L^1(\O)\ .$$

Further verify that $Y_k^\r$ as well as $(Y_k^\r)^{-1}$ are bounded in $L^p(\Omega, \reg) $ independently of $\varrho$. Indeed, for every $\r\in\N$ fixed,
\begin{align}
\lim_{k\to\infty}\int_{\O}|Y^{\r}_{k}|^p\,\md x & =
\int_\O\int_{\reg}|s|^p \,\nu_{x}^{\r}(\md s) \,\md x \nonumber \\
& \le \int_\O\int_{B(0,\r+1)}|s|^p \nu_x(\md s) \,\md
x\le\int_\O\int_{\reg}|s|^p \,\nd \,\md x =\|C\|_{L^1(\O)}
<+\infty,
\label{equiint}
\end{align}
the same calculation could be carried out even $(Y_k^\r)^{-1}$ in place of $Y^{\r}_k$; in this case we need to write $s^{-1}$ in place of $s$.

Applying the diagonalization argument (as $L^1(\O;\C)$ is separable)  we get $\{Y_k\}\subset L^p(\O;\Rm)$ generating $\nu$ and thanks to \eqref{equiint} also \emph{equi-integrable}; the same holds for the inverse.

Moreover, if we defined $\mu$ as the weak* limit of $\mu_\varrho$, then $\mu$ would be generated by $\{Y_k^{-1}\}\subset L^p(\O;\Rm)$ as, due to its definition, 
$$
\lim_{\r\rightarrow\infty}\lim_{k\rightarrow\infty}v((Y^\r_k)^{-1})
=v_{\mu}\mbox{ weakly in }L^1(\O)\ .$$
Also, by applying $\varrho \to \infty$ in \eqref{RelMuRNuR}, it holds that
\be
\int_\reg \hat{v}(s) \nu_x (\d s) = \int_\reg v(s) \mu_x (\d s),
\ee
for all $v \in \C$ and hence, by Lemma \ref{equality} also for all $v \in C^{p,-p}(\reg)$.

\hfill
$\Box$

\bigskip

\section{Proofs of the main results}\label{mains}

\noindent
{\it Proof of Theorem~\ref{THM1}.}
We know from Proposition~\ref{support} that $\nu_x$ is supported on $\reg$ for almost all $x\in\O$. To show that $\nu\in {\cal Y}^{p,-p}(\O;\R^{n\times n})$ we use Proposition~\ref{integrability}. On the other hand, if $\nu\in {\cal Y}^{p,-p}(\O;\R^{n\times n})$ then the existence of a generating sequence is due to Proposition~\ref{explicit}.

It remains to prove relation \eqref{mainpassage}, which we show analogously to \cite[Th.~8.6]{fonseca-leoni}. Let $v:\reg\to\R$ and $g\in L^\infty(\O)$ be as in the theorem.

\noindent For clarity, we divide the proof into 3 steps:

\bigskip

\noindent
{\it Step 1.}
Define $f(x,s):=g(x)v(s)$.  Then $f:\O\times\reg\to\R$ is a normal integrand \cite[Def.~6.27]{fonseca-leoni}.
Suppose first, that $f\ge -M$ for some $M>0$. By \cite[Th.~8.6(i)]{fonseca-leoni} 
\be\label{keyinequality}
\liminf_{k\to\infty}\int_\O f(x,Y_k(x))\,\md x\ge\int_\O\int_\reg f(x,s)\nu_x(\md s)\,\md x\ .
\ee

\bigskip
\noindent
{\it Step 2.}  We use \cite[Th.8.6(i)]{fonseca-leoni} to show that  \eqref{keyinequality} also holds if the negative parts of $f(x,Y_k(x))$, $k\in\N$, form an equi-integrable sequence.   The proof is the same as the proof of \cite[Th.8.6(i)]{fonseca-leoni}. We recall that the negative part of $h:\R\to\R$ is defined as $h^{-}(x):=\max(-h(x),0)$.
 
\bigskip
\noindent
{\it Step 3.}  Here we prove that  if $\{f(x,Y_k(x))\}_{k\in\N}$ is equi-integrable then \eqref{keyinequality} holds as equality.  Namely, if $\{f(x,Y_k(x))\}_{k\in\N}$ is equi-integrable then $\{f^-(x,Y_k(x))\}_{k\in\N}=\{(|f(x,Y_k(x))|-f(x,Y_k(x)))/2\}_{k\in\N}$ is equi-integrable, thus by Step 2:
\be\label{keyinequality0}
\liminf_{k\to\infty}\int_\O f(x,Y_k(x))\,\md x\ge\int_\O\int_\reg f(x,s)\nu_x(\md s)\,\md x\ .
\ee
On the other hand,  $\{-f^{-}(x,Y_k(x))\}_{k\in\N}$ is also equi-integrable, hence
\be\label{keyinequality2}
\liminf_{k\to\infty}\int_\O -f(x,Y_k(x))\,\md x\ge\int_\O\int_\reg -f(x,s)\nu_x(\md s)\,\md x\ .
\ee
Altogether, we have 
\be\label{keyequality}
\lim_{k\to\infty}\int_\O f(x,Y_k(x))\,\md x=\int_\O\int_\reg f(x,s)\nu_x(\md s)\,\md x\ .
\ee

Finally, we show that if $f(x,s)= g(x)v(s)$ for some $g\in L^\infty(\O)$ and $v\in C_{p,-p}(\reg)$ then 
$\{f(x,Y_k(x))\}_{k\in\N}$ is equi-integrable. To see this, we use \cite[Lemma 6.1]{pedregal}. Notice that  $v_0(s):=|v(s)|/(|s|^p+|s^{-1}|^p)\le C$ in $\reg$  for some $C>0$. Moreover, $\lim_{|s|^p+|s^{-1}|^p\to\infty} v_0(s)=0$. Let $(\|Y_k\|_{L^p}^p+\|Y_k^{-1}\|^p_{L^p})\le M$. Take $\varepsilon>0$ and $K>0$ large enough  so that 
$|v_0(s)|<\varepsilon/M$ if $|s|^p+|s^{-1}|^p\ge K/C$. Then for all $k$
\begin{align*}
&\int_{\{x\in\O;\ |v(Y_k(x))|\ge K\}} |v(Y_k(x))|\,\md x \le \int_{\{x\in\O;\ |Y_k(x)|^p+|(Y_k(x))^{-1}|^p \ge K/C\}} |v(Y_k(x))|\,\md x \\ &\quad \leq \int_{\{x \in \Omega; |Y_k(x)|^p+|Y_k^{-1}(x)|^p \geq K/C\}}
|v_0(Y_k(x))|(|Y_k(x)|^p+|Y_k^{-1}(x)|^p) \d x \\ &\quad\leq
\varepsilon/M \int_\Omega |Y_k(x)|^p+|Y_k^{-1}(x)|^p \d x \leq \varepsilon \ .
\end{align*}

\hfill
$\Box$
\bigskip

\noindent
{\it Proof of Theorem~\ref{THM2}.}
It is analogous to the proof of Theorem~\ref{THM1}. Notice that the measure $\nu$  is supported on invertible matrices due to Proposition~\ref{support}. The converse implication follows from Remark~\ref{plus}.
\hfill $\Box$

\bigskip

\noindent{\it Proof of Proposition~\ref{det}.}
By the Mazur lemma det $\nabla y\ge 0$. Suppose that, for contradiction, there existed a set $\omega\subset\O$ of non-zero Lebesgue measure such that det $\nabla y=0$ on $\omega$. We have  by the sequential weak continuity of $W^{1,p}(\O;\R^n)\to L^{p/n}(\O)$: $y\mapsto\det\nabla y$ (\cite{ciarlet}) that
$$
\int_\omega|{\rm det}\ \nabla y_k(x)|\,\md x=\int_\omega{\rm det}\ \nabla y_k(x)\,\md x\to 0\mbox{ as $k\to\infty$}\ , $$
so, it holds for a subsequence (not relabeled) that $0<\det\nabla y_k\to 0$ a.e.~in $\omega$.
By the Fatou lemma, we have
$$
\int_\omega \liminf_{k\to\infty}\frac{\md x}{\det\nabla y_k(x)}\le \liminf_{k\to\infty}\int_\omega\frac{\md x}{\det\nabla y_k(x)}\le C\liminf_{k\to\infty}\int_\omega |(\nabla y_k(x))^{-1}|^n\,\md x\ ,$$
however, the left-hand side tends to $+\infty$. This contradicts the boundedness of $\{(\nabla y_k)^{-1}\}_{k\in\N}$ in $L^p(\O;\R^{n\times n})$ because $p>n$ and $\O$ is bounded. Hence, $\det\nabla y>0$ a.e.~in $\O$. The assertion about the support follows from Proposition~\ref{support}.

\hfill $\Box$

\bigskip

The next two subsections are devoted to the proof of Theorem~\ref{THM3}.

\bigskip

\subsection{Necessary conditions for Young measures to be attained by invertible gradients}
Suppose that we have a bounded  sequence $\{y_k\}\subset W^{1,\infty}(\O;\R^n)$  and  such that $\nabla y_k(x)\in\Rrho$ for some $\varrho>0$ all $k\in\N$ and almost all $x\in\O$. Our aim  is to show that the Young measure generated by $\{\nabla y_k\}$ satisfies  \eqref{supp},\eqref{firstmoment0}, and \eqref{qc0}. In fact, the only point which deserves our attention is the last one, because \eqref{supp} follows easily from Proposition~\ref{Linfty1} and \eqref{firstmoment0} is a well-known description of weak* limits by means of Young measures; cf.~\cite{pedregal}, for instance. 

In the proofs exposed below we will to a large extent follow \cite{k-p1} the main difference compared to this work is that we need to cope carefully  with cut-off functions. Namely, the standard technique of cut-off functions cannot be used since it could destroy the invertibility constraint. 
We will denote by $O(n)$ the set of orthogonal matrices in $\R^{n\times n}$, i.e.,
$O(n):=\{A\in\R^{n\times n};\ A^\top A=AA^{\top}=I\}$ and recall that $\lambda_n(A)$ is the largest singular value of $A\in\R^{n\times n}$, i.e., the largest eigenvalue of $\sqrt{A^\top A}$. We shall heavily rely  on  the following result which can be found in \cite[p.~199 and Remark~2.4]{dacorogna-marcellini}. 

\begin{lemma}\label{daco-marc}
Let $\omega\subset\R^n$ be open and Lipschitz. Let $\varphi\in W^{1,\infty}(\omega;\R^n)$ be such that there is $\vartheta>0$, so that   $0\le \lambda_n(\nabla\varphi)\le 1-\vartheta$ a.e.~ in $\omega$. 
Then there exist   mappings $u\in W^{1,\infty}(\omega;\R^n)$ for which  $\nabla u\in O(n)$ a.e.~in $\omega$ and $u=\varphi$ on $\partial\omega$. Moreover, the set of such  mappings is dense (in the $L^\infty$ norm) in the set $\{\psi:=z+\varphi;\ z\in W_0^{1,\infty}(\omega;\R^n)\ ,\ \lambda_n(\nabla \psi)\le 1-\vartheta \mbox{ a.e.~in } \omega\}$.
\end{lemma}

We shall see in the following that with the aid of this lemma we will be able to construct functions that differ from a particular one only near the boundary. However, this lemma does not allow us to incorporate the bound $\det \nabla y_k > 0$ on minimizing sequences; on sets of arbitrarily small measure we always need to allow also deformations that do not preserve orientation, but still avoiding non-invertibility almost everywhere. Eventually, this technique applies only in the $W^{1,\infty}$-case; for the $W^{1,p}$-case it would be necessary to alter our strategy.

\bigskip

\begin{proposition}\label{boundarycond}
Let $F\in\R^{n\times n}$, $u_F(x):=Fx$ if $x\in\O$,  $y_k\wtostar u_F$ in $W^{1,\infty}(\O;\R^n)$ and let for some $\alpha>0$  $\nabla y_k(x)\in R^{n\times n}_{\alpha}$ for all $k>0$ and almost all $x\in\O$. Then for every $\varepsilon>0$ there is $\{u_k\}\subset W^{1,\infty}(\O;\R^n)$ such that  $\nabla u_k(x)\in R^{n\times n}_{\alpha+\varepsilon}$ for all $k>0$ and almost all $x\in\O$, $u_k-u_F\in W^{1,\infty}_0(\O;\R^n)$ and $|\nabla y_k-\nabla u_k|\to 0$ in measure. In particular, $\{\nabla y_k\}$ and $\{\nabla u_k\}$ generate the same Young measure. 
\end{proposition}

\bigskip

\noindent{\it Proof.} Define for $\ell>0$ sufficiently large  $\O_\ell:=\{x\in\O;\ {\rm dist}(x,\partial\O)\ge 1/\ell\}$ and  smooth cut-off functions $\eta_\ell:\O\to[0,1]$
$$
\eta_\ell(x)=\begin{cases}
1 &\mbox{ if $x\in\O_\ell$}\\
0 &\mbox{ if $x\in\partial\O$}
\end{cases}
$$
such that $|\nabla\eta_\ell|\le C\ell$ for some $C>0$. Define $z_{k\ell}:=\eta_\ell y_k+(1-\eta_\ell)u_F$. Then $z_{k\ell}\in W^{1,\infty}(\O;\R^n)$ and 
$z_{k\ell}=y_k$ in $\O_\ell$ and $z_{k\ell}=u_F$ on $\partial\O$. We see that $\nabla
z_{k\ell}=\eta_\ell\nabla y_k+(1-\eta_\ell)F+ (y_k-u_F)\otimes\nabla \eta_\ell$. Hence, in view of the facts that  $|F|\le\liminf_{k\to\infty}\|\nabla y_k\|_{L^\infty}\le\alpha$ and that $y_k\to u_F$ uniformly in $\bar\O$ we can extract for every $\varepsilon>0$ a (not relabeled) subsequence $k=k(\ell)$ such that 
$$
\|\nabla z_{k(\ell)\ell}\|_{L^\infty} <\alpha +\frac{\varepsilon}{2}\ .$$ 
Consequently, 
$\{z_{k(\ell)\ell}\}$ is uniformly bounded in $W^{1,\infty}(\O;\R^n)$. Moreover,

$$\lambda_n\left(\frac{\nabla z_{k(\ell)\ell}}{\alpha+\varepsilon}\right)\ \le \frac{\|\nabla z_{k(\ell)\ell}\|_{L^\infty}}{\alpha+\varepsilon}\le 1-\frac{\varepsilon}{2(\alpha+\varepsilon)}\ ,$$
where we used the inequality  $\lambda_n(A)\le |A|$ for any $A\in\R^{n\times n}$.
 Denote $\omega_\ell:=\O\setminus\O_\ell$. Then $w_{k(\ell)\ell}:=z_{k(\ell)\ell}|_{\omega_\ell}/(\alpha+\varepsilon)$ is such that $\lambda_n(\nabla w_{k(\ell)\ell})\le 1-\vartheta$ for $\vartheta:=\varepsilon/2(\alpha+\varepsilon)$. 
     We use Lemma~\ref{daco-marc} for $\omega:=\omega_\ell$ and $\varphi:=w_{k(\ell)\ell}$ to  obtain $\phi_{k(\ell)\ell}\in W^{1,\infty}(\omega_\ell;\R^n)$ such that 
   $\phi_{k(\ell)\ell}=w_{k(\ell)\ell}$ on $\partial\omega_\ell$ and $\nabla\phi_{k(\ell)\ell}\in O(n)$. Define 
 $$
u_{k(\ell)\ell}:=\begin{cases}
y_k &\mbox{ if $x\in\O_\ell$}\\
(\alpha+\varepsilon)\phi_{k(\ell)\ell} &\mbox{ if $x\in\O\setminus \O_\ell$.}
\end{cases}
$$  
Notice that $\{u_{k(\ell)\ell}\}_{\ell\in\N}\subset W^{1,\infty}(\O;\R^n)$ and that $u_{k(\ell)\ell}(x)=Fx$ for $x\in\partial\O$. Further, $\nabla u_{k(\ell)\ell}(x)\in R^{n\times n}_{\alpha+\varepsilon}$.    
Moreover, $\nabla u_{k(\ell)\ell}\ne \nabla y_k$ only on sets of vanishing measure, therefore they generate the same Young measure by \cite[Lemma~8.3]{pedregal}.
\hfill
$\Box$ 

\bigskip

\begin{remark}
If $\{u_k\}$ defined in the proof of Proposition~\ref{boundarycond} are homeomorphic and $n=2$ then either  $\det\nabla u_k>0$ or $\det\nabla u_k<0$ in $\O$ for all $k$; cf.~\cite{daneri-pratelli}. The reason is that homeomorphisms in two dimensions are either orientation preserving or reversing. 
\end{remark}

                                                                                                                                                                                                                                                                                                                                  \bigskip
                                                                                                                                                                                                                                                                                                                                  
\begin{lemma}\label{localization}
Let $\nu\in \Grad$. Then for almost all $a\in \O$,  $\mu:=\{\nu_a\}_{x\in\O}\in\Grad$. 
\end{lemma}

\bigskip

{\it Proof.}  We
proceed similarly as  in \cite[Th.~7.2]{pedregal}. Suppose that $\{u_k\}_{k\in\N}\subset W^{1,\infty}(\O;\R^n)$ is a generating sequence for $\nu\in \Grad$ and 
that w$*$-$\lim_{k\to\infty}u_k=u$ in $W^{1,\infty}(\O;\R^n)$. Suppose, moreover, that for some $\varrho > 0$ $\{\nabla u_k(x)\}_{k\in\N} \subset \Rrho$ for a.a. $x \in \O$. 

First we choose $a\in\O$. Define $\bar
V_\ell(y):=\int_{\R^{n\times n}} v^\ell(s)\hat\nu_y(\d s)$ where $\{v^\ell\}_{\ell\in\N}$ is
a dense subset of $\C$. Then we take $a\in\O$, $a\in{\cal
L}_u\cap_{\ell=1}^\infty{\cal
L}_{{V_\ell}}$, where $\mathcal{L}_f$ is the set of Lebesgue points of $f$ in $\O$. The set
of such points has the full Lebesgue  measure.
 Define $u_a:\O\to\R^n$  by  $u_a(x):=\nabla u(a)x$ and denote
$C_a=|\O|^{-1}\int_\O u_a(x)\,\ d x$. Take \be
u^a_{k,j}(x)=j(u_k(a+j^{-1} x)-M_{a,k,j})\ ,\ee where $M_{a,k,j}$
is a constant chosen so that $\int_\O u^a_{k,j}(x)\,\d x=C_a$. Notice that $\nabla u^a_{k,j}(x)= \nabla u_k(a+j^{-1} x)\in\Rrho$ if $j$ is large enough. By
the Poincar\'{e} inequality $\{ u^a_{k,j}\}_{k\in\N,j>0}$ is uniformly
bounded in  $W^{1,\infty}(\O;\R^n)$.

Taking $v\in \C$  and $g\in L^1(\O)$ we have
\begin{eqnarray*}
\int_\O v(\nabla u^a_{k,j}(x))g(x)\,\d x=\int_\O v(\nabla u_k(a+j^{-1} x)g(x)\,\d x
= j^{n}\int_\O v(\nabla
u_k(y))\chi_{a+j^{-1}\O}(y)g\left(\frac{y-a}{j^{-1}}\right)\,\d y\
.
\end{eqnarray*}
Taking now particularly $v^\ell$ for v and using that $\{\nabla y_k\}_{k \in \N}$ generates $\nu$ 
\be\label{limitpass}
\lim_{k\to\infty}\int_\O v^\ell( u^a_{k,j}(x))g(x)\,\d
x&=&j^n\int_\O\bar
V_\ell(y)\chi_{a+j^{-1}\O}(y)g\left(\frac{y-a}{j^{-1}}\right)\,\d y
 \ee
except for all $j \geq j_0 $.
 Passing to the limit for $j\to\infty$ we get by the
Lebesgue dominated convergence theorem 
\begin{eqnarray*}
&&\lim_{j\to\infty}\lim_{k\to\infty}\int_\O v^\ell(\nabla u^a_{k,j}(x))g(x)\,\d x=\lim_{j\to\infty}\int_\O\bar V_\ell(a+j^{-1} x)g(x)\,\d x
=\bar V_\ell (a)\int_\O g(x)\,\d x\\
&=&\int_\O \int_{\R^{n\times n}}v^\ell_0(s)\nu_a(\d s) g(x)\,\d x
=\int_\O \int_{\R^{n\times n}}v^\ell_0(s)\mu_x(\d s) g(x)\,\d x\ .\end{eqnarray*}
The proof is finished by a diagonalization argument.

\hfill $\Box$

\bigskip

We state the next proposition which proves \eqref{qc0}.

 \bigskip

\begin{proposition}\label{proposition:jensen}
Let $\nu\in\Grad$, supp$\,\nu\subset\Rrho$ be  such that for almost all $x\in\O$  $\nabla y(x)= \int_\Rrho s\nu_x(\md s)$, where $y\in W^{1,\infty}(\O;\R^n)$. Then for all $\tilde\varrho\in(\varrho;+\infty]$, almost all $x\in\O$  and all $v\in\mathcal{O}(\tilde\varrho)$ we have 
\be \int_\reg v(s)\nu_x(\md s)\ge Q_{\rm inv}v(\nabla y(x))\ .\ee

\end{proposition}

\bigskip

\noindent{\it Proof.} 
We know from Lemma~\ref{localization} that for almost all $a\in\O$ $\mu=\{\nu_a\}_{x\in\O}\in\Grad$ and that there exits a generating sequence $\{\nabla u_k\}$ such that $\{u_k\}\subset W^{1,\infty}_{\rm inv}(\O;\R^n)$ for $\mu$.  Moreover, $\{u_k\}$  weakly* converges to $x\mapsto\nabla y(a)x$.
Using Proposition \ref {boundarycond} we can, without loss of generality, suppose that for all $k\in\N$ $\nabla u_k\in R^{n\times n}_{\tilde\varrho}$ and $u_k(x)=\nabla y(a)x$ if $x\in\partial\O$ by Lemma~\ref{boundarycond}. Using \eqref{jedna2} for the equality and \eqref{regular-env} the inequality we have

$$
|\O|\int_\reg v(s)\nu_a(\md s) = \lim_{k\to\infty} \int_\O v(\nabla u_k(x))\,\md x \ge |\O|Q_{\rm inv}v(\nabla y(a)) \ .$$
\hfill
$\Box$

\bigskip

\subsection{Sufficient conditions for Young measures to be attained by invertible gradients }
Finally, we show that conditions \eqref{supp},\eqref{firstmoment0}, and \eqref{qc0} are also sufficient for $\nu \in \mathrm{rca}(\reg)$ to be in $\Grad$.
Put 
\be\mathcal{U}^\varrho_A:=\{y\in W^{1,\infty}(\O;\R^n);\ \ y(x)=Ax \  {\rm for  }\ x\in\partial \Omega;\  \max(\|(\nabla y)\|_{L^\infty},\|(\nabla y)^{-1}\|_{L^\infty})\le\varrho\}\ .\ee
Consider for $A\in\R^{n\times n}$ the set 
\be\mathcal{M}^\varrho_A:=\{\overline{\delta_{\nabla y}};\ y\in\mathcal{U}^\varrho_A\}\ ,\ee
where $\overline{\delta_{\nabla y}}\in \rca(\R^{n\times n})$ is defined as
\be
\la \overline{\delta_{\nabla y}}, v\ra:=|\O|^{-1}\int_\O v(\nabla y(x))\,\md x\ .\ee

We have the following  lemma.

\begin{lemma}\label{convexity}
 Let $A\in\R^{n\times n}$  If $\varrho>|A|$ then the set $\mathcal{M}^\varrho_A$ is nonempty and  convex.
\end{lemma}

\bigskip

\noindent {\it Proof.} First we show that $\mathcal{M}^\varrho_A$ is non-empty. This is clear when $A$ is invertible, otherwise  we use Lemma \ref{daco-marc}. Take $A\in\R^{n\times n}$ and $\varrho>|A|$.   We recall that $\lambda_n(A/\varrho)$ is the largest singular value of $A/\varrho$. We have  $\lambda_n(A/\varrho)\le |A|/\varrho=1-(\varrho-|A|)/\varrho$. Apply  Lemma~\ref{daco-marc} with $\varphi(x):=Ax/\varrho$, $x\in\O$ and $\vartheta:=(\varrho-|A|)/\varrho$  to get  $u\in W^{1,\infty}(\O;\R^n)$ such that $\nabla u \in O(n)$ a.e.~in $\O$ and $u(x)=Ax/\varrho$ if $x\in\partial\O$. Therefore, $y:=\varrho u\in\mathcal{U}^\varrho_A $. Consequently, $\mathcal{M}^\varrho_A\ne\emptyset$. 

The rest of proof  is analogous to  the proof of \cite[Lemma~8.5]{pedregal}.  We take $y_1,y_2\in \mathcal{U}^\varrho_A$ and for a given $\lambda\in(0,1)$ we find a subset $D\subset\O$ such that $|D|=\lambda|\O|$. There are two countable families  of subsets of $D$ and $\O\setminus D$ of the form 
$$
\{a_i+\epsilon_i\O;\ a_i\in D,\ \epsilon_i>0,\ a_i+\epsilon_i\O\subset D\}$$
and 
$$
\{b_i+\epsilon_i\O;\ b_i\in \O\setminus D,\ \rho_i>0,\ b_i+\rho_i\O\subset\O \setminus D\}\ $$
such that 
$$
D=\cup_{i}(a_i+\epsilon_i\O)\cup N_0 \ ,$$ 

$$\O\setminus D=\cup_{i}(b_i+\rho_i\O)\cup N_1  \ ,$$ 
 
where the Lebesgue measure of $N_0$ and $N_1$ is zero. 
We define 

\be
y(x):=
\begin{cases}
\epsilon_iy_1\left(\frac{x-a_i}{\epsilon_i}\right)+Aa_i & \mbox{ if $x\in a_i+\epsilon_i\O$, }\\
 \rho_iy_2\left(\frac{x-b_i}{\rho_i}\right)+Ab_i & \mbox{ if $x\in b_i+\rho_i\O$, }\\
 Ax &\mbox{ otherwise.}
\end{cases}
\ee

Then 
\be
\nabla y(x)=
\begin{cases}
\nabla y_1\left(\frac{x-a_i}{\epsilon_i}\right)& \mbox{ if $x\in a_i+\epsilon_i\O$, }\\
 \nabla y_2\left(\frac{x-b_i}{\rho_i}\right) & \mbox{ if $x\in b_i+\rho_i\O$. }\\
\end{cases}
\ee

In particular, $y\in{\mathcal U}_A$ and 
$$\overline{\delta_{\nabla y}}=\lambda\overline{\delta_{\nabla y_1}}+(1-\lambda)\overline{\delta_{\nabla y_2}}\ .$$

Notice that  $\|(\nabla y_i)^j\|_{L^\infty(\O;\R^n)}\le \varrho$, $i=1,2$, $j=-1,1$, and therefore $\|(\nabla y)^j\|_{L^\infty(\O;\R^n)}\le \varrho$, as well.
\hfill
$\Box$

\bigskip

\begin{remark}\label{finiteness}
It follows from the proof of the previous Lemma that if $v\in C(\reg)$  then $Q_{\rm inv} v<+\infty$; cf.~\eqref{regular-env}.

\end{remark}

\bigskip

\begin{lemma}\label{homogenization}
Let  $\{u_k\}_{k\in\N}\subset \mathcal{U}_A$ from \eqref{UF} be a bounded sequence. Let
$\nu\in \Grad$ be generated by $\{\nabla u_k\}$. Then
there is a bounded sequence $\{w_k\}\subset \mathcal{U}_A$  such that
 $\{\nabla w_k\}_{k\in\N}$ generates
$\bar\nu\in \Grad$,
and  for any $v\in \C$ and almost all $x\in\O$
\be\label{homog} \int_{\R^{n\times n}}
 v(s)\overline{\nu}_x(\md s)=
\frac{1}{|\O|}\int_{\O}\int_{\R^{n\times n}} v(s)\nu_x(\md s)\,\md x\ . \ee

\end{lemma}

\bigskip

\noindent{\it Sketch of proof.} We follow the proof of \cite[Th.~7.1]{pedregal}. The
family
$$
{\cal A}=\left\{x\in a+\epsilon\bar\O\subset\O\ ;\ a\in\O\ ,\
\epsilon\le j^{-1}\right\}\ $$ is a  covering of $\O$. There
exists a countable collection $\{x\in
a_{ij}+\epsilon_{ij}\bar\O\}$, $\epsilon_{ij}\le 1/j$ of pairwise
disjoint sets  and
$$\O=\bigcup_i\{x\in  a_{ij}+\epsilon_{ij}\bar\O\}\bigcup N_j\ ,\ |N_j|=0\ .$$
We see that $\sum_{i}\epsilon_{ij}^n=|\O|/|\O|=1$. We now take for
$u_A(x)=Ax$, $x\in\O$, the following sequence of mappings
$$
w_k(x) := \left\{ \begin{array}{ll}
           \epsilon_{ik}u_k\left(\frac{x-a_{ik}}{\epsilon_{ik}}\right) +u_A(a_{ik})  & \mbox{ if $x\in a_{ik}+\epsilon_{ik}\O$}\\
                            u_A(x) & \mbox{otherwise}\ .
                            \end{array}
                   \right.
$$
Therefore, $w_k=u_Y$ on $\partial\O$ and for a.a. $x\in\O$
$$
\nabla w_k(x)=\nabla
u_k\left(\frac{x-a_{ik}}{\epsilon_{ik}}\right)\ .$$ 
Hence, the Poincar\'e inequality yields the bound on $\{w_k\}$ in
$W^{1,\infty}(\O;\R^m)$ and we even see that $\{w_k\}$ is bounded in $\mathcal{U}_A$.  See the proof of \cite[Th.~7.1]{pedregal} to verify that $\{\nabla w_k\}$ generates $\bar\nu$.

\hfill $\Box$

\begin{proposition}\label{homocase}
Let $\mu$ be a probability measure supported on a compact set  $K\subset \R^{n\times n}_\alpha$ for some $\alpha>0$  and let  $ A:=\int_K s\mu(\md s)$.  Let $\varrho>\alpha$ and let 
\be\label{jensen} 
Q_{\rm inv}v(A)\le \int_{K} v(s)\mu(\md s)\ ,\ee
for all $v\in \mathcal{O}(\varrho)$.
Then $\mu\in \Grad$ and it is generated by gradients of mappings from $\mathcal{U}^\varrho_A$.
\end{proposition}

\bigskip

{\it Proof.}
The proof standardly uses the Hahn-Banach theorem and Lemma~\ref{homogenization} and it is  similar to  \cite[Proposition~8.17]{pedregal}. First, notice that $|A|\le \alpha<\varrho $. Then since $\mathcal{M}^\varrho_A$ is non-empty and convex due to Lemma~\ref{convexity}, we can by the Hahn-Banach theorem find a continuous linear functional $T:\rca(\Rrho)\to \R$ such that $T(\nu)\ge 0$ on $\bar{\mathcal{M}}^\varrho_A$ and $T(\nu) < 0$ otherwise. Using the Riesz representation theorem we therefore find a $\tilde v\in C(\Rrho)$ such that 
$$
0\le \la T(\nu),\tilde v\ra=\int_\Rrho \tilde v(s)\nu(\md s) =|\O|^{-1}\int_\O \tilde v(\nabla y(x))\,\md x\ , $$
for all $\nu \in \mathcal{M}^\varrho_A$ and hence all $y\in\mathcal{U}_A^\varrho$.
If 
$$
v(s):=\begin{cases}
\tilde v(s) &\mbox{ if $s\in\Rrho$,}\\
+\infty &\mbox{otherwise}
\end{cases}
$$
then $v\in\mathcal{O}(\varrho)$ and $0 \leq Q_{\rm inv} v(A)=\inf_{\mathcal{U}_A}|\O|^{-1}\int_\O \tilde v(\nabla y(x))\,\md x$.
By \eqref{jensen}, $0\le \int_{\Rrho} v(s)\mu(\md s)$. Thus, $\mu\in \overline{\mathcal{M}^\varrho_A}$ (the weak* closure). As $C(\Rrho)$ is separable, weak* topology of bounded sets in $\rca(\Rrho)$ is   metrizable. 
Hence, there is a sequence $\{u_k\}_{k\in\N}\subset \mathcal{U}_A^\varrho$  such that for all $v\in C(\Rrho)$ (and all $v\in\mathcal{O}(\varrho)$)
$$
\lim_{k\to\infty}\int_\O v(\nabla u_k(x))\,\md x= |\O|\int_{\Rrho} v(s)\mu(\md s)\ ,$$
and $\{u_k\}$ is bounded in $W^{1,\infty}(\O;\R^{n\times n})$ due to the Poincar\'{e} inequality.  As $u_k(x)=Ax$ for $x\in\partial\O$ we use the homogenization procedure 
from Lemma~\ref{homogenization} to show that $\mu$ is the homogeneous Young measure generated by $\{\nabla u_k\}$.

\hfill
$\Box$

\bigskip

We will need the following auxiliary result.

\begin{lemma}\label{auxiliary} (see \cite[Lemma~6.1]{k-p1})
Let $\O\subset\R^n$ be an open domain  with $|\partial\O|=0$ and
let $N\subset\O$ be of the zero Lebesgue measure. For
$r_k:\O\setminus N\to (0,+\infty)$ and $\{f_k\}_{k\in\N}\subset
L^1(\O)$ there exists a set of points $\{a_{ik}\}\subset
\O\setminus N$ and positive numbers $\{\epsilon_{ik}\}$,
$\epsilon_{ik}\le r_k(a_{ik})$ such that
$\{a_{ik}+\epsilon_{ik}\bar\O\}$ are pairwise disjoint for each
$k\in\N$, $\bar\O=\cup_i \{a_{ik}+\epsilon_{ik}\bar\O\}\cup N_k$
with $|N_k|=0$ and for any $j\in\N$ 
$$
\lim_{k\to\infty}\sum_i
f_j(a_{ik})|\epsilon_{ik}\O|= \int_\O
f_j(x)\,\md x\ .$$

\end{lemma}

\bigskip

The next proposition forms the sufficiency part of Theorem~\ref{THM3}.

\bigskip

\begin{proposition}
Let $\O\subset\R^n$, let $\nu=\{\nu_x\}_{x\in\O}$ be a family of probability measures on $\R^{n\times n}$.  Suppose that for some $\varrho>0$ and  for almost all $x\in\O$  supp $\nu_x\subset\Rrho$. Let, moreover, the following 
two conditions hold: \\
\be\label{firstmoment1}
 \exists\ u\in W^{1,\infty}(\O;\R^n)\ :\  \nabla u(x)=\int_{\reg} s \nu_x(\d s)\  \ ,\ee
for all  $\tilde\varrho>\varrho$, for a.a. $x\in\O$ and all  $v\in \mathcal{O}(\tilde\varrho)$  the following  inequality
is valid \be\label{qc2} Q_{\rm inv}v(\nabla u(x))\le\int_{\reg} v(s)\nu_x(\md s)\ .
 \ee Then
$\nu\in\Grad$.

\end{proposition}

\bigskip

\noindent {\it Proof.} Some parts of the proof follow 
\cite[Proof of Th.~6.1]{k-p1}. 
We are looking for a sequence $\{u_k\}_{k\in\N}\subset
W^ {1,\infty}(\O;\R^n)$ satisfying
$$
\lim_{k\to\infty} \int_\O v(\nabla u_k(x))g(x)\,\md
x=\int_{\O} \int_{\R^{n\times
n}}v(s)\nu_x(\md s)g(x)\,\md x\  $$ for
all $g\in \Gamma$ and any $v\in S$, where
$\Gamma$ and $S$ are countable dense subsets of $C(\bar\O)$ and
$\C$, respectively.

First of all notice that, as $u\in W^ {1,\infty}(\O;\R^n)$ it is differentiable in $\O$ outside a set of measure zero called $N$, we may find for every $a\in\O\setminus N$ and every $k > 0$ a $r_k(a)>0$ for any $0<\epsilon < r_k(a)$ we have
\be\label{derivative}
 \frac{1}{\epsilon}| u(a+\epsilon y)-u(a)-\epsilon \nabla u(a)y|\le \frac1k\ .
\ee
Furthermore, due to the continuity of $g$ we choose  $r_k(a)>0$ smaller if necessary to assure that for any  $0<\epsilon < r_k(a)$ 
\begin{equation}
\label{g-est}
\int_{a+\epsilon\O}g(x)\,\md x = g(a)\epsilon + \frac{1}{k}. 
\end{equation} 
From Lemma~\ref{auxiliary} we can find find $a_{ik}\in\O\setminus N$, $\epsilon_{ik}\le  r_k(a_{ik})$ ($r_k(a_{ik})$ are defined above by equation \eqref{derivative}) such that for all $v\in S$ and all $g\in \Gamma$
\be\label{79} \lim_{k\to\infty}\sum_i \bar
V(a_{ik})g(a_{ik}) |\epsilon_{ik} \O|= \int_\O \bar
V(x)g(x)\,\md x\ ,\ee
where
$$\bar V(x):=\int_{\reg} v(s)\nu_x(\md s)\ .$$
In view of  Lemma~\ref{homocase}, we can assume  that $\{\nu_{a_{ik}}\}_{x\in\O}$ is a homogeneous gradient Young measure living in $\Grad$ and we call $\{\nabla u^{ik}_j\}_{j\in\N}$ its generating sequence. We know that we can consider  $\{u^{ik}_j\}_j\subset \mathcal{U}^{\tilde\varrho}_{\nabla u(a_{ik})}$ for arbitrary $\tilde\varrho>\varrho$. So, it yields

\be\label{imp14} \lim_{j\to\infty} \int_\O v(\nabla
u_j^{ik}(x))g(x)\,\md x=\bar V(a_{ik})\int_\O g(x)\,\md x\ \ee
and, in addition, 
\be\label{wl} {\rm
w}^*-\lim_{j\to\infty}u^{ik}_j=\nabla u(a_{ik})x \mbox{ in $W^ {1,\infty}(\O;\R^n)$}
\ee

Let $\O_\ell:=\{x\in\O;\ {\rm dist}(x,\partial\O)\ge \ell^{-1}\}$.
We define a sequence of smooth cut-off functions
$\{\eta_\ell\}_{\ell\in\N}$ such that
$$
\eta_\ell(x): = \left\{ \begin{array}{ll}
          0   & \mbox{ in $\O_\ell$}\ ,\\
                            1 & \mbox{on $\partial\O$}\
                            \end{array}
                   \right.
$$
and $|\nabla\eta_\ell|\le C\ell$ for some $C>0$.

Further, take  a sequence $\{u_k^\ell\}_{k,\ell\in\N}\subset
W^ {1,\infty}(\O;\R^n)$ defined  by
$$
u_k^\ell(x): = \left\{ \begin{array}{ll}
           \left[u(a_{ik})+\epsilon_{ik}u^{ik}_j\left(\frac{x-a_{ik}}{\epsilon_{ik}}\right)\right]\left(1-\eta_\ell\left(\frac{x-a_{ik}}{\epsilon_{ik}}\right)\right)\\
+ u(x)\eta_\ell\left(\frac{x-a_{ik}}{\epsilon_{ik}}\right)   & \mbox{ if $x\in a_{ik}+\epsilon_{ik}\O$,}\\
                            u(x) & \mbox{otherwise}\ ,
                            \end{array}
                   \right.
$$
where $j=j(i,k,\ell)$ will be chosen later. Note that for every
$k$  we have $u^\ell_k-u\in W_0^ {1,\infty}(\O,\R^n)$.

We calculate for $x\in a_{ik}+\epsilon_{ik}\O$
\begin{eqnarray}\label{grad}
\nabla u_k^\ell(x)&=&\nabla u^{ik}_j\left(\frac{x-a_{ik}}{\epsilon_{ik}}\right)\left(1-\eta_\ell\left(\frac{x-a_{ik}}{\epsilon_{ik}}\right)\right)\nonumber\\
&+& \nabla u(x)\eta_\ell\left(\frac{x-a_{ik}}{\epsilon_{ik}}\right)\nonumber\\
&+& \frac{1}{\epsilon_{ik}}\left[u(x)-u(a_{ik})-\epsilon_{ik}\nabla u(a_{ik})\left(\frac{x-a_{ik}}{\epsilon_{ik}}\right)\right]\otimes\nabla\eta_\ell\left(\frac{x-a_{ik}}{\epsilon_{ik}}\right)\nonumber\\
&+&\left[\nabla u(a_{ik})\left(\frac{x-a_{ik}}{\epsilon_{ik}}\right)-u^{ik}_j\left(\frac{x-a_{ik}}{\epsilon_{ik}}\right)\right]\otimes\nabla\eta_\ell\left(\frac{x-a_{ik}}{\epsilon_{ik}}\right).
\end{eqnarray}
Notice that moduli of all four terms can be made together  uniformly bounded by $\tilde\varrho>\varrho$. Namely, notice that the sum of  the first two terms is $\le\varrho$ and the other two terms can be made arbitrarily small if $k$ is sufficiently large compared to $\ell$ by exploiting \eqref{derivative} and the strong convergence in $L^\infty(a_{ik} + \epsilon_{ik}\O; \R^n)$ of $u^{ik}_j(x)$ to $\nabla u(a_{ik})x$.  

Take the set $(a_{ik}+\varepsilon_{ik}\Omega)\setminus
(a_{ik}+\varepsilon_{ik}\Omega_\ell)$ and solve the inclusion  $\nabla \tilde u^\ell_k\in O(n)$ with the boundary conditions $\tilde u^\ell_k= u^\ell_k/\tilde\varrho$ if $x\in \partial ((a_{ik}+\varepsilon_{ik}\Omega_k)\setminus(a_{ik}+\varepsilon_{ik}\Omega_\ell))$. This inclusion has a solution due to Lemma~\ref{daco-marc}. Set

 $$
z_k^\ell(x): = \left\{ \begin{array}{ll}
           u_k^\ell(x)  & \mbox{ if $x\in a_{ik}+\epsilon_{ik}\O_\ell$,}\\
           \tilde u_k^\ell(x) &\mbox { if $x\in (a_{ik}+\epsilon_{ik}\O)\setminus (a_{ik}+\epsilon_{ik}\O_\ell)$,}\\
                            u(x) & \mbox{otherwise}\ .
                            \end{array}
                   \right.
$$
Observe, that the Lebesgue measure of  the set $\{x\in\O;\ \nabla (u^\ell_k(x)-z^\ell_k(x))\ne 0\}$  vanishes as $\ell\to\infty$. Further, $\{z_k^\ell\}_{k,\ell}\subset W^{1,\infty}(\O;\R^n)$ is a bounded sequence as well as $\{\nabla z_k^\ell\}_{k,\ell}^{-1}\subset L^\infty(\O;\R^{n\times n})$.

Let us fix $k,i,\ell$ (with $k$ sufficiently large such that the $|\nabla z_k^\ell|$ is uniformly bounded by $\tilde{\rho}$) and consider sets $\{E_k\}_{k\in\N}$, $E_k\subset E_{k+1}$ and $\Gamma\times S=\cup_{k}E_k$ . We can eventually enlarge each $j=j(i,k,\ell)$ so  that additionally for
any $(g,v_0)\in E_k$ \be\label{line} \left|\epsilon_{ik}^n\int_\O
g(a_{ik}+\epsilon_{ik}y)v(\nabla u_j^{ik}(y))\,\md y -\bar
V(a_{ik})\int_{a_{ik}+\epsilon_{ik}\O}g(x)\,\md
x\right|\le\frac{1}{2^ik}\ .\ee

We have for some $C>0$ 
\begin{eqnarray*}
&&\int_\O g(x)v(\nabla u^\ell_k(x))\,\md x =\sum_i \epsilon_{ik}^n\int_\O g(a_{ik}+\epsilon_{ik}y)v(\nabla u_j^{ik}(y))\,\md y +\frac{C}{\ell}\ .
\end{eqnarray*}
Here we used the smallness of $|\O\setminus\O_\ell| $ as $\ell\to\infty$ and boundedness of $g$ and $v$.

Consequently, in view of \eqref{79}, \eqref{g-est} and \eqref{line}  for all
$(g,v)\in\Gamma\times S$
$$ \lim_{\ell\to\infty}\lim_{k\to \infty} \int_\O g(x)v(\nabla u^\ell_k(x))\,\md x=\int_\O\int_{\R^{\times n}} v(s)\nu_x(\md s)g(x)\, \md x\ .$$
Hence, we can pick up a subsequence $\{\nabla u^{\ell}_{k(\ell)}\}_\ell$ generating $\nu$. The measure $\nu$ is also generated by $\{\nabla z^{\ell}_{k(\ell)}\}_\ell$ because the difference of both sequences vanishes in measure.
Finally, we see from the construction  that  $\{z^\ell_{k(\ell)}\}_\ell$ can be chosen to have the same boundary conditions
as $u$.

\hfill$\Box$

 \bigskip

{\bf Acknowledgment:}  BB was supported by the grants P201/10/0357 (GA \v{C}R), 41110 (GAUK \v{C}R), and the research plan AV0Z20760514 (\v CR). MK also thanks for support through the grant P201/12/0671 (GA \v{C}R). GP wishes to thank GA \v{C}R for the support through the project 
 P105/11/0411.

\end{sloppypar}
\end{document}